# RUIN PROBABILITIES AND OVERSHOOTS FOR GENERAL LÉVY INSURANCE RISK PROCESSES[1]

By Claudia Klüppelberg, Andreas E. Kyprianou and Ross A. Maller

*Munich University of Technology, University of Utrecht and Australian National University*

We formulate the insurance risk process in a general Lévy process setting, and give general theorems for the ruin probability and the asymptotic distribution of the overshoot of the process above a high level, when the process drifts to $-\infty$ a.s. and the positive tail of the Lévy measure, or of the ladder height measure, is subexponential or, more generally, convolution equivalent. Results of Asmussen and Klüppelberg [*Stochastic Process. Appl.* **64** (1996) 103–125] and Bertoin and Doney [*Adv. in Appl. Probab.* **28** (1996) 207–226] for ruin probabilities and the overshoot in random walk and compound Poisson models are shown to have analogues in the general setup. The identities we derive open the way to further investigation of general renewal-type properties of Lévy processes.

**1. Introduction.** Various recent studies of insurance risk processes and associated random walks and Lévy processes have paid particular attention to the heavy-tailed case, when downward jumps of the process—claims— may be very large. Such models are now thought to be quite realistic, especially in view of a recent tendency to large-claim events in the insurance industry.

To give some intuition for the much more general framework of this paper, we briefly recall the classical insurance risk model, where all quantities are explicit. In the classical model, the claims arriving within the interval $(0, t]$,

Received February 2003; revised September 2003.
[1]Supported in part by ARC Grant DP0210572.

$t > 0$, are modelled as a compound Poisson process, yielding the risk process

$$R_t = u + \gamma t - \sum_{i=1}^{N_t} Y_i, \qquad t \geq 0, \tag{1.1}$$

where $u$ is the initial risk reserve and $\gamma > 0$ is the premium rate (as usual, we set $\sum_{i=1}^{0} a_i = 0$). Denote by $F$ the claim size distribution function (d.f.), that is, the d.f. of the independent and identically distributed (i.i.d.) almost surely (a.s.) positive random variables (r.v.'s) $Y_i$, assumed to have finite mean $\mu > 0$. Let $\lambda > 0$ be the intensity of the Poisson process, assume $\gamma > \lambda\mu$, and let $\rho = \lambda\mu/\gamma < 1$. The probability of ultimate ruin is then

$$\psi(u) = P(R_t < 0 \text{ for some } t > 0)$$

$$= P\left(\sum_{i=1}^{N_t} Y_i - \gamma t > u \text{ for some } t > 0\right)$$

$$= P\left(\sum_{i=1}^{n} (Y_i - \gamma T_i) > u \text{ for some } n \in \mathbb{N}\right) \tag{1.2}$$

$$= (1 - \rho) \sum_{n=1}^{\infty} \rho^n \overline{F_I^{*n}}(u). \tag{1.3}$$

We have used the following notation and facts:

(i) In this model, ruin can occur only when a claim occurs. This, jointly with the fact that the interarrival times $\{T_i : i \in \mathbb{N}\}$ of a Poisson process are i.i.d. exponential r.v.'s, leads to (1.2).

(ii) Equation (1.3) follows from a ladder height analysis; in this classical case, the *integrated tail distribution*

$$F_I(x) := \frac{1}{\mu} \int_0^x \overline{F}(y)\,dy, \qquad x \geq 0, \tag{1.4}$$

is the d.f. associated with the increasing ladder height process of the process $X_t = \sum_{i=1}^{N_t} Y_i - \gamma t$, $t \geq 0$, and $\overline{F_I^{*n}}$ is the tail of its $n$-fold convolution.

(iii) The condition $\rho < 1$ guarantees that the process $X$ has negative drift.

(iv) The infinite sum in (1.3) constitutes a renewal measure, which is defective with killing rate $\rho$.

All this standard theory can be found in various textbooks, for example, [20] and [2], to mention just the classic and the most recent one.

In analyzing (1.3), two regimes can be recognized. The first is called the *Cramér case*, when there exists a $\nu > 0$ satisfying $\rho \int_0^\infty e^{\nu u} F_I(dx) = 1$. The defect $\rho$ in the renewal function (1.3) can then be removed by an exponential tilting, and, using Smith's key renewal lemma (see, e.g., [16], Section 1.2), the ruin probability can be shown to decrease exponentially fast, in fact, proportional to $e^{-\nu u}$, as $u \to \infty$. This result has been extended to a Lévy process setting by Bertoin and Doney [6].

If such a "Lundberg coefficient" $\nu$ does not exist, as is the case for subexponential and other "convolution equivalent" distributions (see Section 3), estimates of the ruin probability have been derived by Embrechts, Goldie and Veraverbeke [15], Embrechts and Veraverbeke [18] and Veraverbeke [35]; see [16], Section 1.4.

It is also of prime interest to understand the way ruin happens. This question has been addressed by Asmussen [1] for the Cramér case, and more recently by Asmussen and Klüppelberg [3] for the subexponential case. They describe the sample path behavior of the process along paths leading to ruin via various kinds of conditional limit theorems. As expected, the Cramér case and the non-Cramér case are qualitatively quite different; see, for example, [2] and [16], Section 8.3.

Our aim is to investigate the non-Cramér case in a general Lévy process setting, which clearly reveals the roles of the various assumptions. Our Lévy process $X$ will start at 0 and be assumed to drift to $-\infty$ a.s., but otherwise is quite general. Upward movement of $X$ represents "claim payments," and the drift to $-\infty$ reflects the fact that "premium income" should outweigh claims. "Ruin" will then correspond to passage of $X$ above a specified high level, $u$, say. In this scenario, heavy-tailedness of the positive side of the distribution of upward jumps models the occurrence of large, possibly ruinous, claims, and has previously been studied in connection with the assumption of a finite mean for the process. But in general we do not want to restrict the process in this way. A higher rate of decrease of the process to $-\infty$ is more desirable from the insurer's point of view, while allowing a heavier tail for the positive part is in keeping with the possibility of even more extreme events, which indeed are observed in recent insurance data.

This leads to the idea of considering processes for which the only assumption is of a drift to $-\infty$ a.s., possibly at a linear rate, as is the case when the mean is finite and negative, but possibly much faster. This kind of analysis is aided by results going back to [19] which allow us to quantify such behavior, as is done, for example, via easily verified conditions for drift to $\pm\infty$ given in [11]. We also make essential use of important fluctuation identities given in [4] and [36]. Our results can thus be seen as adding to an understanding of renewal and fluctuation properties of Lévy processes which drift to $-\infty$, with application to passage time and overshoot behavior at high levels.

The paper is organized as follows. In the next two sections we introduce some basic notation, definitions and results for later use in the study. These consist of some renewal-theoretic aspects of Lévy processes (drifting to $-\infty$) in Section 2, together with definitions and properties of subexponential and related classes of distributions in Section 3. In Section 4 we present our main results, which concern the asymptotic analyses of first passage times and the ruin probability, asymptotic conditional overshoot distributions and some ladder height and ladder time considerations. Section 5 establishes some useful asymptotic relations between the Lévy measures of $X$ and its ladder height process, while Section 6 offers some examples of the results presented in Section 4. Proofs of the main results are given in Section 7.

**2. Some renewal theory for Lévy processes.** Let us suppose that $X = \{X_t : t \geq 0\}$ is a general Lévy process with law $P$ and Lévy measure $\Pi_X$. That is to say, $X$ is a Markov process with paths that are right continuous with left limits such that the increments are stationary and independent and whose characteristic function at each time $t$ is given by the Lévy–Khinchine representation

$$E(e^{i\theta X_t}) = e^{-t\Psi(\theta)}, \qquad \theta \in \mathbb{R},$$

where

$$(2.1) \quad \Psi(\theta) = i\theta a + \sigma^2\theta^2/2 + \int_{(-\infty,+\infty)} (1 - e^{i\theta x} + i\theta x \mathbf{1}_{\{|x|<1\}}) \Pi_X(dx).$$

We have $a \in \mathbb{R}$, $\sigma^2 \geq 0$ and $\Pi_X$ is a measure supported on $\mathbb{R}$ with $\Pi_X(\{0\}) = 0$ and $\int_\mathbb{R} (x^2 \wedge 1) \Pi_X(dx) < \infty$ ([4], page 13, and [32], Chapter 2). The natural filtration generated by $X$ is assumed to satisfy the usual assumptions of right continuity and completeness.

Throughout we impose three essential restrictions:

(i) $X_0 = 0$ and the process drifts to $-\infty$: $\lim_{t \to \infty} X_t = -\infty$ a.s;
(ii) $\Pi_X\{(0, \infty)\} > 0$, so the process is not spectrally negative;
(iii) we consider the non-Cramér case [see (4.3) and Proposition 5.1].

Further discussion of these points is given below. Otherwise, the only requirement will be on the asymptotic tail behavior (convolution equivalence, see Definition 3.2) which we assume for the right tail of $\Pi_X$.

The following are standard tools of fluctuation theory for Lévy processes; see, for example, [4], Chapter VI.

DEFINITION 2.1.

*Supremum.* Let $\overline{X} = \{\overline{X}_t = \sup_{s \in [0,t]} X_s : t \geq 0\}$ be the process of the last supremum.

*Local time and inverse local time.* Let $L = \{L_t : t \geq 0\}$ denote the *local time* in the time period $[0, t]$ that $\overline{X} - X$ spends at zero. Then $L^{-1} = \{L_t^{-1} : t \geq 0\}$ is the *inverse local time* such that

$$L_t^{-1} = \inf\{s \geq 0 : L_s > t\}.$$

We shall also understand
$$L_{t-}^{-1} = \inf\{s \geq 0 : L_s \geq t\}.$$
In both cases and in the following text, we take the infimum of the empty set as $\infty$. Note that the previous two inverse local times are both stopping times with respect to the natural filtration of $X$. Since $X$ drifts to $-\infty$, it follows that, with probability 1, $L_\infty < \infty$ and hence there exists a $t > 0$ such that $L_t^{-1} = \infty$, again with probability 1.

*Increasing and decreasing ladder height processes.* The process $H$ defined by $\{H_t = X_{L_t^{-1}} : t \geq 0\}$ is the *increasing ladder height process*, that is to say, the process of new maxima indexed by local time at the maximum. We call $L^{-1}$ the (*upwards*) ladder time process. The processes $L^{-1}$ and $H$ are both defective subordinators. It is understood that $H_t = \infty$ when $L_t^{-1} = \infty$.

We shall define $\hat{H} = \{\hat{H}_t : t \geq 0\}$ to be the *decreasing ladder height process* in an analogous way. Note that this means that $\hat{H}$ is a process which is negative valued (this is unconventional, as the usual definition of decreasing ladder height process would correspond to $-\hat{H}$ here).

*Bivariate ladder process.* Given the event $\{0 \leq t < L_\infty\}$, the joint process $(L^{-1}, H)$ behaves on $[0, t)$ like a bivariate subordinator which is independent of $L_\infty$. Also there exists a constant $q > 0$ such that $L_\infty \stackrel{d}{=} \mathbf{e}_q$, where $\mathbf{e}_q$ is an exponential variable with mean $1/q$; compare [4], Lemma VI.2. Throughout the paper we shall distinguish between the nondefective processes, denoted by $\mathcal{L}$ (with $\mathcal{L}_\infty = \infty$), $\mathcal{L}^{-1}$ and $\mathcal{H}$, and their defective versions $L$, $L^{-1}$ and $H$. The corresponding nondefective bivariate ladder process is then $(\mathcal{L}^{-1}, \mathcal{H})$. It is a bivariate subordinator, independent of $\mathbf{e}_q$, with the property

(2.2) $$\{(L_t^{-1}, H_t) : t < L_\infty\} \stackrel{\text{Law}}{=} \{(\mathcal{L}_t^{-1}, \mathcal{H}_t) : t < \mathbf{e}_q\}.$$

Note that, by contrast, the decreasing ladder height process is not defective in this sense because we have assumed that $X$ drifts to $-\infty$.

DEFINITION 2.2 (Lévy measures and their tails). In addition to the measure $\Pi_X$, we shall denote by $\Pi_\mathcal{H}$ and $\Pi_{\hat{H}}$ the Lévy measures of $\mathcal{H}$ and $\hat{H}$, with supports in $(0, \infty)$ and $(-\infty, 0)$, respectively. Further, for $u > 0$,
$$\overline{\Pi}_X^+(u) = \Pi_X\{(u, \infty)\},$$
$$\overline{\Pi}_X^-(u) = \Pi_X\{(-\infty, -u)\},$$
$$\overline{\Pi}_X = \overline{\Pi}_X^+(u) + \overline{\Pi}_X^-(u)$$
represent the positive, negative and combined tails of $\Pi_X$. We use analogous notation for the tails of $\Pi_\mathcal{H}$ and $\Pi_{\hat{H}}$.

In our applications, the first passage time

$$\tau(u) = \inf\{t \geq 0 : X_t > u\}, \qquad u > 0,$$

corresponds to ruin occurring at level $u$, and major objects of interest are the probability that this occurs in a finite time, and the behavior of this probability as the reserve level $u$ is increased to $\infty$. Following [3] and [15], a natural way to proceed is by placing subexponential or, more generally, "convolution equivalence" assumptions (see Section 3) on $\overline{\Pi}_{\mathcal{H}}$ or on $\overline{\Pi}_X^+$. We are then able to follow in outline the program of [3], finding the limiting conditional distribution as $u \to \infty$ of the overshoot $X_{\tau(u)} - u$ above level $u$ (when it occurs), and of further quantities in our general setup. This gives quite a clear picture of how and when first passage over a high level happens for general Lévy processes.

The following development is essentially based on appropriate sections of [4] and [36], but adapted and extended in part for our requirements.

DEFINITION 2.3 (Ladder height renewal measure). We define the renewal measure, $V$, of the defective process $H$ in the usual way. Its connection to the nondefective process $\mathcal{H}$ with exponential killing time is as follows:

$$(2.3) \quad dV(y) = \int_0^\infty dt \cdot P(H_t \in dy) = \int_0^\infty dt \cdot e^{-qt} P(\mathcal{H}_t \in dy), \qquad y \geq 0.$$

We shall also be interested in the renewal measure, $\hat{V}$, of $\hat{H}$, the downward ladder height process, satisfying

$$d\hat{V}(y) = \int_0^\infty dt \cdot P(\hat{H}_t \in dy), \qquad y \leq 0.$$

The next theorem gives an identity from which we can calculate the distributions of $X_{\tau(u)}$, $L^{-1}_{L_{\tau(u)^-}}$ and $X_{L^{-1}_{L_{\tau(u)^-}}}$. Although notationally rather complicated, the latter two objects are nothing more than the time corresponding to the ladder time prior to the first passage time (i.e., to the ruin time), and the position of this ladder height, respectively.

THEOREM 2.4.  *Fix $u > 0$. Suppose that $f$, $g$ and $h$ are bounded, positive and measurable, and that $g(u) = 0$. Define*

$$dV^h(y) = \int_0^\infty dt \cdot e^{-qt} \int_{[0,\infty)} h(\phi) P(\mathcal{H}_{t^-} \in dy, \mathcal{L}_t^{-1} \in d\phi), \qquad y \geq 0.$$

*Then*

$$(2.4) \quad \begin{aligned} &E\Big(f\big(X_{L^{-1}_{L_{\tau(u)^-}}}\big) g(X_{\tau(u)}) h\big(L^{-1}_{L_{\tau(u)^-}}\big); \tau(u) < \infty\Big) \\ &\qquad = \int_{(0,u]} dV^h(y) f(y) \int_{(u-y,\infty)} g(y+s) \Pi_{\mathcal{H}}(ds). \end{aligned}$$

PROOF. Define $T(u) = \inf\{t \geq 0 : H_t > u\}$ and recall that $X$ experiences first passage at $\tau(u)$ if and only if $H$ experiences first passage at $T(u)$. The quantity $X_{L^{-1}_{L_{\tau(u)-}}}$ can alternatively be written as $H_{T(u)-}$. On $\{T(u) < \infty\}$, $H$ is a subordinator and $L_\infty$ has an Exponential($q$) distribution. Start from the left-hand side of the statement of the theorem, and decompose according to $\{T(u) = t\}$ to get

$$E\Big(f\big(X_{L^{-1}_{T(u)-}}\big)g(H_{T(u)})h(L^{-1}_{T(u)-}); T(u) < L_\infty\Big)$$
$$= E\sum_{0 < t < L_\infty}(f(H_{t-})g(H_{t-} + \Delta H_t)h(L^{-1}_{t-}); H_{t-} \leq u < H_{t-} + \Delta H_t)$$
$$= E\int_0^\infty dy \cdot qe^{-qy}\sum_{0 < t < y}(f(\mathcal{H}_{t-})g(\mathcal{H}_{t-} + \Delta\mathcal{H}_t)h(\mathcal{L}^{-1}_{t-});$$
$$\mathcal{H}_{t-} \leq u < \mathcal{H}_{t-} + \Delta\mathcal{H}_t)$$
$$= E\sum_{t>0}e^{-qt}(f(\mathcal{H}_{t-})g(\mathcal{H}_{t-} + \Delta\mathcal{H}_t)h(\mathcal{L}^{-1}_{t-}); \mathcal{H}_{t-} \leq u < \mathcal{H}_{t-} + \Delta\mathcal{H}_t).$$

Use the compensation formula for the Poisson point process $\{\triangle\mathcal{H}_t : t \geq 0\}$ ([4], page 7) to get that the last expression is equal to

$$\int_0^\infty dt \cdot e^{-qt}E\Big(f(\mathcal{H}_{t-})h(\mathcal{L}^{-1}_{t-})\mathbf{1}_{\{\mathcal{H}_{t-} \leq u\}}\int_{(0,\infty)}\Pi_\mathcal{H}(ds)g(\mathcal{H}_{t-} + s)\mathbf{1}_{\{\mathcal{H}_{t-}+s>u\}}\Big)$$
$$= \int_{(0,u]}\Big\{\int_0^\infty dt \cdot e^{-qt}\int_{[0,\infty)}h(\phi)P(\mathcal{H}_{t-} \in dy, \mathcal{L}^{-1}_{t-} \in d\phi)\Big\}$$
$$\times f(y)\int_{(u-y,\infty)}\Pi_\mathcal{H}(ds)g(y+s)$$
$$= \int_{(0,u]}dV^h(y)f(y)\int_{(u-y,\infty)}\Pi_\mathcal{H}(ds)g(y+s). \qquad \square$$

The proof of Theorem 2.4 is similar to calculations appearing in Proposition III.2, page 76, [4] (see also [37]). The seemingly curious condition $g(u) = 0$ functions as a way of excluding from the calculation the fact that there is possibly an atom at $u$ in the distribution of $X_{\tau(u)}$ which is a result of crossing $u$ continuously or "creeping upwards" (see Remark 2.8).

The next result, giving a formula for the ruin probability, is the continuous time version of the Pollacek–Khinchine formula (see [4], page 172, and [6], page 364).

PROPOSITION 2.5. $P(\tau(u) < \infty) = qV(u,\infty) := q\overline{V}(u), \ u > 0.$

DEFINITION 2.6 (Wiener–Hopf factors). The Wiener–Hopf factorization theorem (see, e.g., [4], page 166), together with the downward drift assumption on our Lévy process, tells us that we can write, for some constant $k > 0$,

$$\begin{aligned}(2.5) \quad k\Psi(\theta) &= -k\log Ee^{i\theta X_1} = [-\log Ee^{i\theta H_1}][-\log Ee^{i\theta \hat{H}_1}] \\ &= \kappa(\theta) \times \hat{\kappa}(\theta), \qquad \theta \in \mathbb{R}.\end{aligned}$$

The constant $k$ is determined by the choice of normalization of the local time $L$. We may and will assume without loss of generality that $k = 1$. A different value of $k$ would simply modify the choices of $\mathcal{L}$, $\mathcal{H}$ and $q$. We have, for $\nu > 0$ and some $c \geq 0$, $\hat{c} \geq 0$,

$$(2.6) \qquad \kappa(i\nu) = -\log Ee^{-\nu H_1} = \int_{(0,\infty)} (1 - e^{-\nu y})\Pi_{\mathcal{H}}(dy) + \nu c + q$$

and, recalling that $\hat{H}$ is negative,

$$(2.7) \qquad \hat{\kappa}(-i\nu) = -\log Ee^{\nu \hat{H}_1} = \int_{(-\infty,0)} (1 - e^{\nu y})\Pi_{\hat{H}}(dy) + \nu \hat{c}.$$

The factors $\kappa(\cdot)$ and $\hat{\kappa}(\cdot)$ are the Lévy–Khinchine exponents of $H$ and $-\hat{H}$, which are subordinators, and accordingly the integrals in the definitions of $\kappa$ and $\hat{\kappa}$ converge. The nonnegative constants $c$ and $\hat{c}$ are the drift coefficients of these subordinators and $q$ is the same killing rate that appears in the definition of $\mathcal{H}$ [see (2.2)]. The convention that $e^{i\theta H_1} = 0 = e^{-\nu H_1}$ when $H_1 = \infty$ is used in (2.5) and (2.6).

REMARK 2.7. Since $Ee^{-\nu H_1} = e^{-q}Ee^{-\nu \mathcal{H}_1}$ for $\nu > 0$, (2.6) implies

$$(2.8) \qquad -\log Ee^{-\nu \mathcal{H}_1} = \int_{(0,\infty)} (1 - e^{-\nu y})\Pi_{\mathcal{H}}(dy) + \nu c, \qquad \nu > 0,$$

and, as a consequence of (2.5), (2.6) and (2.8), we have, for $\nu > 0$,

$$q - \log Ee^{-\nu \mathcal{H}_1} = \frac{\Psi(i\nu)}{\hat{\kappa}(i\nu)},$$

and hence

$$(2.9) \qquad q = \lim_{\nu \downarrow 0} \frac{\Psi(i\nu)}{\hat{\kappa}(i\nu)}.$$

The limit in (2.9) exists, and can be easily calculated, for example, when $X_1$ has finite mean, in which case $\Psi$ and $\hat{\kappa}$ are differentiable at 0; see the examples in Section 6.

REMARK 2.8 (Creeping). $X$ is said to *creep upward* if $P(X_{\tau(u)} = u, \tau(u) < \infty) > 0$ for some (hence every) $u > 0$; equivalently, if the $c$ defined in (2.6) is positive ([4], pages 174 and 175). $X$ creeps downward if $-X$ creeps upward; equivalently, if the $\hat{c}$ defined in (2.7) is positive. Suppose $c > 0$. Then we have

$$P(X_{\tau(u)} = u, \tau(u) < \infty) = P(T'(u) < L_\infty) = E(e^{-qT'(u)}), \qquad u > 0,$$

where $T'(u) = \inf\{t \geq 0 : H_t = u\}$. A similar proof as in Theorem 5, page 79, of [4], applied to the defective subordinator $H$, then shows that the derivative $dV(u)/du$ exists and is continuous and positive on $(0, \infty)$, and that

$$(2.10) \qquad P(X_{\tau(u)} = u, \tau(u) < \infty) = c \frac{dV}{du}(u) =: cV'(u), \qquad u > 0.$$

When $c = 0$, $V'$ is not defined, but the next corollary (to Theorem 2.4), which lists the main formulae that we will use, shows that we do not need it then.

COROLLARY 2.9. *We have the following four convolution identities for* $u > 0$:

(i) $P(X_{\tau(u)} - u > x, \tau(u) < \infty) = \int_{(0,u]} dV(y)\overline{\Pi}_\mathcal{H}(u + x - y)$;

(ii) $P(\tau(u) < \infty) = \int_{(0,u)} dV(y)\overline{\Pi}_\mathcal{H}(u - y) + cV'(u)$, *with the convention that the term containing* $V'(u)$ *is absent when* $V'$ *is not defined, that is, when* $c = 0$;

(iii) $P(X_{\tau(u)} > u, L^{-1}_{L_{\tau(u)}-} > \psi, \tau(u) < \infty) = \int_{(0,u)} V(dy; \psi)\overline{\Pi}_\mathcal{H}(u-y)$, *where* $V(dy; \psi) = \int_0^\infty dt \cdot e^{-qt} P(\mathcal{H}_t \in dy, \mathcal{L}_t^{-1} > \psi)$;

(iv) $P(X_{L^{-1}_{L_{\tau(u)}-}} > \phi, \tau(u) < \infty) = \int_{(\phi,u)} V(dy)\overline{\Pi}_\mathcal{H}(u-y) + cV'(u), \phi \in [0, u)$, *again with the convention that the term containing* $V'(u)$ *is absent when* $V'$ *is not defined*.

PROOF. (i) Just choose $f = h = 1$ and $g = \mathbf{1}_{\{\cdot > x + u\}}$ in Theorem 2.4.

(ii) Multiply each side of the equation in (ii) by $e^{-\nu u}$, with $\nu > 0$, and integrate over $u \in [0, \infty)$, making use of Proposition 2.5 and the identities

$$(2.11) \qquad \int_{[0,\infty)} e^{-\nu y} V(dy) = \frac{1}{q - \log E e^{-\nu \mathcal{H}_1}}$$

[obtained by integrating (2.3)], and

$$(2.12) \quad \begin{aligned} &\nu c + \nu \int_{(0,\infty)} e^{-\nu y} \overline{\Pi}_\mathcal{H}(y)\, dy \\ &= \nu c + \int_{(0,\infty)} (1 - e^{-\nu y})\Pi_\mathcal{H}(dy) = -\log E e^{-\nu \mathcal{H}_1} \end{aligned}$$

[from (2.8)] to see that (ii) holds as stated. [Note that by taking the limit as $x$ tends to zero in (i), and combining the result with (ii), we recover (2.10).]

(iii) Choosing $f = 1$, $g = \mathbf{1}_{\{\cdot > x+u\}}$ and $h = \mathbf{1}_{\{\cdot > \psi\}}$ in (2.4), and taking the limit as $x$ tends to zero, gives (iii).

(iv) Choosing $f = \mathbf{1}_{\{\cdot > \phi\}}$, $g = \mathbf{1}_{\{\cdot > x+u\}}$ and $h = 1$ in (2.4), then letting $x$ tend to zero, gives an expression for $P(X_{L^{-1}_{L_{\tau(u)^-}}} > \phi, \tau(u) < \infty, X_{\tau(u)} > u)$. Since $\{X_{L^{-1}_{L_{\tau(u)^-}}} = u\}$ on $\{X_{\tau(u)} = u\}$ almost surely, by adding on $cV'$ we have (iv). □

Further convolution identities that will be of use can be found in Proposition 3.3 of [36].

THEOREM 2.10 (Vigon [36]). *We have, for $u \in (0, \infty)$:*

(i) $\overline{\Pi}^+_X(u) = \int_{(u,\infty)} \overline{\Pi}_{\hat{H}}(u-y)\, d\Pi_{\mathcal{H}}(y) + \hat{c}\Pi'_{\mathcal{H}}(u)$, *where $\Pi'_{\mathcal{H}}$ is the density of $\Pi_{\mathcal{H}}$, which exists if and only if $\hat{c}$, the drift coefficient of $-\hat{H}$, is positive;*

(ii) $\overline{\Pi}_{\mathcal{H}}(u) = -\int_{(-\infty,0)} \overline{\Pi}^+_X(u-y)\, d\hat{V}(y)$.

REMARK 2.11. Note that by our convention $\hat{V}(y)$ is positive and nonincreasing on $y \in (-\infty, 0)$, with $\hat{V}(0) = 0$. $X$ drifts to $-\infty$ a.s. in our analysis, so we can and will exclude the case when $X$ is a subordinator. This means that $\hat{H}$, $\hat{V}$ and $\hat{\kappa}$ are not identically zero.

We say that $X$ is *spectrally positive* if $\Pi_X\{(-\infty, 0)\} = 0$. We then have $\hat{H}_t = -t$ and hence $\hat{V}(dy) = -dy$ and $\hat{c} = 1$, and the expressions in Vigon's theorem simplify considerably. In particular, (i) and (ii) both say that

$$(2.13) \qquad \overline{\Pi}_{\mathcal{H}}(u) = \int_u^\infty \overline{\Pi}^+_X(y)\, dy = \int_u^\infty \overline{\Pi}_X(y)\, dy, \qquad u > 0$$

[further implying that the integral in (2.13) is finite, and thus also that $E|X_1|$ is finite]. See [4], Chapter VII, for other useful results concerning spectrally one-sided processes.

**3. Convolution equivalence and infinite divisibility.** Each infinitely divisible d.f. generates a Lévy process in the sense that it may serve as the d.f. of $X_1$. For the most part we shall restrict ourselves to those infinitely divisible d.f.'s which belong to one of the following classes.

DEFINITION 3.1 (Class $\mathcal{L}^{(\alpha)}$). Take a parameter $\alpha \geq 0$. We shall say that a d.f. $G$ on $[0, \infty)$ with tail $\overline{G} := 1 - G$ belongs to class $\mathcal{L}^{(\alpha)}$ if $\overline{G}(x) > 0$ for each $x \geq 0$ and

$$(3.1) \quad \lim_{u \to \infty} \frac{\overline{G}(u-x)}{\overline{G}(u)} = e^{\alpha x} \qquad \text{for each } x \in \mathbb{R}, \text{ if } G \text{ is nonlattice;}$$

(3.2) $$\lim_{n\to\infty} \frac{\overline{G}(n-1)}{\overline{G}(n)} = e^\alpha \qquad \text{if } G \text{ is lattice (then assumed of span 1)}.$$

(There should be no confusion of the class $\mathcal{L}^{(\alpha)}$ with our notation $\mathcal{L}_t$ for the local time.)

DEFINITION 3.2 (Convolution equivalence and class $\mathcal{S}^{(\alpha)}$). With $*$ denoting convolution, $G$ is said to be *convolution equivalent* if $G \in \mathcal{L}^{(\alpha)}$ for some $\alpha \geq 0$, and if in addition, for some $M < \infty$, we have

(3.3) $$\lim_{u\to\infty} \frac{\overline{G^{*2}}(u)}{\overline{G}(u)} = 2M,$$

where $\overline{G^{*2}}(u) = 1 - G^{*2}(u)$. We say that $G$ belongs to $\mathcal{S}^{(\alpha)}$. The class $\mathcal{S}^{(0)}$ is called the *subexponential distributions*. The parameter $\alpha$ is referred to as the *index* of the class $\mathcal{S}^{(\alpha)}$ (or $\mathcal{L}^{(\alpha)}$). We will often write $\overline{G} \in \mathcal{L}^{(\alpha)}$ rather than $G \in \mathcal{L}^{(\alpha)}$, and similarly for $\mathcal{S}^{(\alpha)}$.

A number of useful properties flow from these definitions. The limit relation (3.1) holds locally uniformly. In [14] it is shown that, when $G \in \mathcal{S}^{(\alpha)}$, then any d.f. $F$ which is tail equivalent to $G$ [i.e., $\overline{F}(x) \sim \overline{G}(x)$ as $x \to \infty$, equivalently $\lim_{x\to\infty} \overline{F}(x)/\overline{G}(x) = 1$] is also in $\mathcal{S}^{(\alpha)}$. The tail of any (Lévy or other) measure, finite and nonzero on $(x_0, \infty)$ for some $x_0 > 0$, can be renormalized to be the tail of a d.f., and, by extension, then is said to be in $\mathcal{L}^{(\alpha)}$ or $\mathcal{S}^{(\alpha)}$ if the appropriate conditions in Definitions 3.1 or 3.2 are satisfied. For these results and others, see, for example, [14, 15], and their references.

We follow Bertoin and Doney [7] in (3.1) and (3.2). They drew attention to the need, when $\alpha > 0$, to distinguish the lattice and nonlattice cases; under (3.1), the geometric distribution, for example, would not be in $\mathcal{L}^{(\alpha)}$. For $\alpha = 0$, no distinction is necessary. Having noted this distinction for $\alpha > 0$, we will confine our proofs to the nonlattice case by considering (3.1) to be the defining property of $\mathcal{L}^{(\alpha)}$.

DEFINITION 3.3 (Moment generating function). For a finite d.f. $G$ on $[0, \infty)$, the *moment generating function* is defined (for all $a \in \mathbb{R}$ such that the following integral is finite) as

$$\delta_a(G) = \int_{[0,\infty)} e^{au} G(du).$$

Of course, $\delta_0(G) < \infty$. When $G \in \mathcal{S}^{(\alpha)}$ for an $\alpha > 0$, Fatou's lemma applied to (3.3), using (3.1), shows that $\delta_\alpha(G) < \infty$. Furthermore, the constant $M$ in (3.3) must then equal $\delta_\alpha(G)$ (cf. [9, 10, 31]). Moreover, $\delta_{\alpha+\varepsilon}(G) = \infty$ for all $\varepsilon > 0$. For the class $\mathcal{S}^{(0)}$ of subexponential d.f.'s, the latter property means that the moment generating function does not exist for any $\varepsilon > 0$—these distributions are "heavy-tailed" in this sense. Typical examples are Pareto, heavy-tailed Weibull and lognormal d.f.'s. Distributions with regularly varying tails are in this class. Note that while the Exponential($\alpha$) distribution itself is in $\mathcal{L}^{(\alpha)}$ (for the same index $\alpha$), it is not in $\mathcal{S}^{(\beta)}$ for any index $\beta \geq 0$; the convolution of two Exponential($\alpha$) distributions is a Gamma$(2, \alpha)$ distribution for which (3.3) does not hold. Distributions in the class $\mathcal{S}^{(\alpha)}$ for $\alpha > 0$ are, however, "near to exponential" in the sense that their tails are only slightly modified exponential; see [27]. The slight modification, however, results in a moment generating function which is finite for argument $\alpha$, as observed above. An important class of d.f.'s which are convolution equivalent or subexponential for some values of the parameters is the generalized inverse Gaussian distributions, having densities

$$f(x) = \left(\frac{b}{a}\right)^{d/2} (2K_d(\sqrt{ab}))^{-1} x^{d-1} \exp\left(-\frac{1}{2}(ax^{-1} + bx)\right), \qquad x > 0,$$

where $K_d$ is the modified Bessel function of the third kind with index $d$. The following parameter sets are possible: $\{a \geq 0, b > 0\}$ for $d \geq 0$; $\{a > 0, b > 0\}$ for $d = 0$; $\{a > 0, b \geq 0\}$ for $d < 0$. (For $a = 0$ or $b = 0$, the respective limits are to be taken in the norming constants.) For this distribution, $F \in \mathcal{L}^{(b/2)}$ for each $b \geq 0$, and, when $d < 0$, $F \in \mathcal{S}^{(b/2)}$ for each $b \geq 0$; see [13] and [28].

Extending (3.3), when $G \in \mathcal{S}^{(\alpha)}$ for an $\alpha \geq 0$, it is in fact true that, for all $k \in \mathbb{N}$,

$$(3.4) \qquad \lim_{u \to \infty} \frac{\overline{G^{*k}}(u)}{\overline{G}(u)} = k\delta_\alpha^{k-1}(G).$$

Also, the following uniform bound due to Kesten holds: for each $\varepsilon > 0$, there is a $K(\varepsilon)$ such that, uniformly in $u > 0$ for each $k \in \mathbb{N}$,

$$(3.5) \qquad \frac{\overline{G^{*k}}(u)}{\overline{G}(u)} \leq K(\varepsilon)(\delta_\alpha(G) + \varepsilon)^k.$$

An important property of $\mathcal{S}^{(\alpha)}$ relates these classes to infinitely divisible distributions, and hence to Lévy processes.

PROPOSITION 3.4. *Fix an $\alpha \geq 0$. If $G$ is infinitely divisible with Lévy measure $\Pi_G(\cdot) \neq 0$, whose tail is $\overline{\Pi}_G(u) = \Pi_G\{(u, \infty)\}$, $u > 0$, then the following are equivalent:*

(3.6)
(i) $\overline{G} \in \mathcal{S}^{(\alpha)}$;
(ii) $\overline{\Pi}_G \in \mathcal{S}^{(\alpha)}$;
(iii) $\overline{\Pi}_G \in \mathcal{L}^{(\alpha)}$ *and* $\lim_{u \to \infty} \dfrac{\overline{G}(u)}{\overline{\Pi}_G(u)} = \delta_\alpha(G).$

For a proof of Proposition 3.4 in the case $\alpha = 0$, see Embrechts, Goldie and Veraverbeke [15]; they restrict themselves to distributions on $[0, \infty)$, while Pakes [30] gives the result for distributions on $(-\infty, \infty)$, and for $\alpha \geq 0$. For more detailed information on the classes $\mathcal{S}^{(\alpha)}$, and in particular on the subexponential class, we refer to [16] and the review paper [24].

The next lemma applies Proposition 3.4 to get some basic asymptotic relations for the tail of the ladder height process $\mathcal{H}_t$ and for the ruin probability.

LEMMA 3.5. *Fix an $\alpha \geq 0$. Suppose $P(\mathcal{H}_1 > u) \in \mathcal{S}^{(\alpha)}$, or, equivalently, by Proposition 3.4, $\overline{\Pi}_{\mathcal{H}} \in \mathcal{S}^{(\alpha)}$. Then for each $t > 0$,*

$$(3.7) \quad P(\mathcal{H}_t > u) \sim t\delta_\alpha^t(\mathcal{H})\overline{\Pi}_{\mathcal{H}}(u) \sim t\delta_\alpha^{t-1}(\mathcal{H})P(\mathcal{H}_1 > u), \qquad u \to \infty,$$

*and hence, by tail equivalence, $P(\mathcal{H}_t > u) \in \mathcal{S}^{(\alpha)}$ for each $t > 0$. Suppose further that $e^{-q}\delta_\alpha(\mathcal{H}) < 1$. Then*

$$(3.8) \qquad \lim_{u \to \infty} \frac{P(\tau(u) < \infty)}{\overline{\Pi}_{\mathcal{H}}(u)} = \frac{q}{(q - \log \delta_\alpha(\mathcal{H}))^2} = q\delta_\alpha^2(V).$$

[Here and throughout, we write $\delta_\alpha(\mathcal{H})$ for $\delta_\alpha(\mathcal{H}_1)$.]

PROOF. Apply (3.6) to the infinitely divisible r.v. $\mathcal{H}_t$ with Lévy measure $\Pi_{\mathcal{H}_t}(\cdot) = t\Pi_{\mathcal{H}}(\cdot)$ to get, for each $t > 0$,

$$P(\mathcal{H}_t > u) \sim \delta_\alpha(\mathcal{H}_t)\overline{\Pi}_{\mathcal{H}_t}(u) = t\delta_\alpha^t(\mathcal{H})\overline{\Pi}_{\mathcal{H}}(u), \qquad u \to \infty,$$

then apply (3.6) again to complete (3.7). Next, use the fact that $P(\mathcal{H}_t > u)$ does not decrease in $t$ (for each $u > 0$) and the discrete uniform bound (3.5) to see that for each $\varepsilon > 0$, there is a $K(\varepsilon)$ such that, for all $t > 0$ and $u > 0$,

$$(3.9) \quad \begin{aligned} P(\mathcal{H}_t > u) &\leq P(\mathcal{H}_{\lfloor t \rfloor + 1} > u) \\ &\leq K(\varepsilon)(\delta_\alpha(\mathcal{H}) + \varepsilon)^{\lfloor t \rfloor + 1} P(\mathcal{H}_1 > u), \qquad u > 0. \end{aligned}$$

Proposition 2.5 gives

$$\frac{P(\tau(u) < \infty)}{\overline{\Pi}_{\mathcal{H}}(u)} = \frac{q\overline{V}(u)}{\overline{\Pi}_{\mathcal{H}}(u)} = \frac{q}{\overline{\Pi}_{\mathcal{H}}(u)}\int_0^\infty e^{-qt}P(\mathcal{H}_t > u)\,dt, \qquad u > 0.$$

Then (3.7) and the uniform bound (3.9), together with dominated convergence, and assuming that $e^{-q}\delta_\alpha(\mathcal{H}) < 1$, give

$$\lim_{u \to \infty} \frac{P(\tau(u) < \infty)}{\overline{\Pi}_{\mathcal{H}}(u)} = q\int_0^\infty e^{-qt}t\delta_\alpha^t(\mathcal{H})\,dt = \frac{q}{(q - \log \delta_\alpha(\mathcal{H}))^2}.$$

The final equality in (3.8) follows from (2.11), as we can put $\nu = -\alpha$ when $e^{-q}\delta_\alpha(\mathcal{H}) < 1$. $\square$

**4. Main results.** Throughout the entire paper we assume

(4.1) $\qquad X_0 = 0, \qquad \lim_{t\to\infty} X_t = -\infty \text{ a.s.}, \qquad \Pi_X\{(0,\infty)\} > 0.$

[The spectrally negative case, when $\Pi_X\{(0,\infty)\} = 0$, is easily dealt with separately in our context; see Remark 4.6.]

Our main assumption throughout this section will be

(4.2) $\qquad\qquad\qquad\qquad \overline{\Pi}_{\mathcal{H}} \in \mathcal{S}^{(\alpha)},$

for a specified $\alpha \geq 0$.

For the specified $\alpha$, the *non-Cramér condition*,

(4.3) $\qquad\qquad\qquad\qquad e^{-q}\delta_\alpha(\mathcal{H}) < 1,$

will also be assumed in our main results. This condition has force only when $\alpha > 0$; for $\alpha = 0$, condition (4.3) is automatically satisfied when (4.1) holds, since $q > 0$ then.

We start with the asymptotics of the first passage time $\tau(u)$ in Theorem 4.1, which extends Lemma 3.5 by showing that (3.8) can only hold if the ruin probability is in $\mathcal{S}^{(\alpha)}$.

THEOREM 4.1 (Limiting first passage time, $\alpha \geq 0$). *Fix an $\alpha \geq 0$ and assume (4.1)–(4.3) hold. Then, as $u \to \infty$,*

(4.4) $\begin{aligned} P(\tau(u) < \infty) &\sim \frac{q}{(q - \log \delta_\alpha(\mathcal{H}))^2} \overline{\Pi}_{\mathcal{H}}(u) \\ &\sim \frac{q}{(q - \log \delta_\alpha(\mathcal{H}))^2 \delta_\alpha(\mathcal{H})} P(\mathcal{H}_1 > u), \end{aligned}$

*and thus $\overline{C}(u) := P(\tau(u) < \infty)$, $u > 0$, is in $\mathcal{S}^{(\alpha)}$. Conversely, suppose that (4.1) holds and $\overline{C}(\cdot)$ is in $\mathcal{S}^{(\alpha)}$. Then (4.2) and (4.3), and hence (4.4), hold.*

To be practically useful, we need to replace the quantities depending on the ladder variables in Theorem 4.1 (and similarly in our other results) with quantities defined as far as possible in terms of the marginal distributions of $X$ or, better, in terms of $\Pi_X$. Section 5 is devoted to results like this so we defer discussion until then.

THEOREM 4.2 (Overshoot, local time at ruin, last ladder height before ruin, $\alpha \geq 0$). *Fix an $\alpha \geq 0$ and assume (4.1)–(4.3) hold. Then:*

(i) *for all $x > 0$,*

(4.5) $\qquad\qquad \lim_{u\to\infty} P(X_{\tau(u)} - u > x \mid \tau(u) < \infty) = \overline{G}(x),$

*where $\overline{G}$ is the tail of a (possibly improper) distribution function:*

(4.6) $\quad \overline{G}(x) = \frac{e^{-\alpha x}}{q}\left(q - \log \delta_\alpha(\mathcal{H}) + \int_{(x,\infty)} (e^{\alpha y} - e^{\alpha x})\Pi_{\mathcal{H}}(dy)\right), \qquad x \geq 0;$

(ii) *for all $t \geq 0$,*

(4.7)
$$\lim_{u \to \infty} P(L_{\tau(u)} > t \mid \tau(u) < \infty) \\ = e^{-(q - \log \delta_\alpha(\mathcal{H}))t}(1 + t(q - \log \delta_\alpha(\mathcal{H})) \log \delta_\alpha(\mathcal{H})/q);$$

(iii) *for all $\phi \geq 0$,*

(4.8) $$\lim_{u \to \infty} P\Big(X_{L^{-1}_{L_{\tau(u)^-}}} \leq \phi \mid \tau(u) < \infty\Big) = \frac{(q - \log \delta_\alpha(\mathcal{H}))^2}{q} \left( \int_{(0,\phi]} e^{\alpha y} V(dy) \right).$$

REMARK 4.3. (i) In the last result, when $\alpha = 0$, the limiting distribution is proper. This follows by virtue of the fact that

$$V(\infty) = \int_0^\infty e^{-qt} dt = 1/q.$$

On the other hand, when $\alpha > 0$, the limiting distribution is improper, having mass at infinity

$$1 - \delta_\alpha(V) \frac{(q - \log \delta_\alpha(\mathcal{H}))^2}{q} = 1 - \frac{1}{q \delta_\alpha(V)} > 0.$$

(ii) When $\alpha > 0$, we can let $x \to 0+$ in (4.6) to see that $\overline{G}(0^+) = 1 - \alpha c/q$; thus we can also conclude that the asymptotic conditional probability of creeping over the barrier $u$, as $u \to \infty$, is equal to $\alpha c/q$. When ruin occurs, the process has positive probability of crossing the boundary by creeping or jumping.

(iii) When $\alpha = 0$, the distribution $G$ in Theorem 4.2 is degenerate, placing all its mass at infinity. Ruin thus occurs asymptotically only by a jump.

For the case $\alpha = 0$, we have the following sharper result:

THEOREM 4.4 (Sharper limiting overshoot distribution, $\alpha = 0$). *Suppose that (4.1) holds and $\overline{\Pi}_\mathcal{H} \in \mathcal{S}^{(0)}$. Then, for all $x > 0$,*

(4.9) $$\lim_{u \to \infty} \left| P(X_{\tau(u)} - u > x \mid \tau(u) < \infty) - \frac{\overline{\Pi}_\mathcal{H}(u + x)}{\overline{\Pi}_\mathcal{H}(u)} \right| = 0,$$

*and the convergence is uniform in $x \geq \eta$ for each $\eta > 0$.*

The remaining result in this section concerns the last ladder time before ruin. For this, we only show tightness:

PROPOSITION 4.5 (Last ladder time before ruin, $\alpha = 0$). *Assume (4.1) and $\overline{\Pi}_\mathcal{H} \in \mathcal{S}^{(0)}$. Then*

$$\lim_{\phi \to \infty} \limsup_{u \to \infty} P(L^{-1}_{L_{\tau(u)^-}} > \phi \mid \tau(u) < \infty) = 0.$$

REMARK 4.6 (Spectrally negative case). In this case, $\Pi_X\{(0,\infty)\} = 0$, and there are no upward jumps, so we have $X_{\tau(u)} = u$ on $\tau(u) < \infty$, for all $u > 0$, $X$ creeps up, and the overshoot is a.s. zero at all levels. The ladder height process $\mathcal{H}_t$ is simply the unit drift $t$ ([4], page 191). The passage time $\tau(u)$ has Laplace transform

$$E(e^{-\lambda\tau(u)}; \tau(u) < \infty) = e^{-u\Phi(\lambda)},$$

where $\Phi$ is the right inverse function to $-\Psi(-i\lambda)$ ([4], page 189). Thus the ruin probability is $P(\tau(u) < \infty) = e^{-au}$, where $a > 0$ satisfies $\Psi(-ia) = 0$. In the classical risk model, this setup is taken to describe annuities in life insurance ([25], page 9).

**5. Relations between $\Pi_X$, $\Pi_\mathcal{H}$ and $q$.** In this section we give some useful connections between the m.g.f.'s and the Lévy measures of $X$ and $\mathcal{H}$, and related quantities.

PROPOSITION 5.1 [Criteria for (4.3)]. *Assume* (4.1). *For any $\nu > 0$, the following equivalences are true:*

$$\begin{aligned}
E(e^{\nu X_1}) \ &(\text{is finite and}) \ < 1 \\
&\iff \ e^{-q}\delta_\nu(\mathcal{H}) < 1 \iff \delta_\nu(V) < \infty \\
&\iff \log\delta_\nu(\mathcal{H}) = \nu c + \int_{[0,\infty)} (e^{\nu y} - 1)\Pi_\mathcal{H}(dy) < q \\
&\iff \nu a - \sigma^2\nu^2/2 - \int_{(-\infty,\infty)} (e^{\nu x} - 1 - \nu x \mathbf{1}_{\{|x|<1\}})\Pi_X(dx) > 0,
\end{aligned}$$

(5.1)

*and if any of the conditions holds then*

$$(5.2) \qquad \frac{1}{\delta_\nu(V)} = q - \log\delta_\nu(\mathcal{H}) = \frac{-\log Ee^{\nu X_1}}{-\log Ee^{\nu \hat{H}_1}}.$$

REMARK 5.2. In the case $\alpha > 0$, Proposition 5.1 shows that our results in Section 4 apply to the class of Lévy processes for which $Ee^{\alpha X_1} < 1$. By contrast, suppose there is a $\nu_0 > 0$ such that $Ee^{\nu_0 X_1} = 1$. This forces $X$ to drift to $-\infty$ a.s., and, without further assumptions, Bertoin ([4], page 183) and Bertoin and Doney [6] then prove Cramér's estimate: $P(\tau(u) < \infty) \sim Ce^{-\nu_0 u}$, as $u \to \infty$, where $C < \infty$, and $C > 0$ if and only if the Lévy process $X^\#$ with exponent $\Psi^\#(\lambda) = \Psi(\lambda - i\nu_0)$ has $E|X_1^\#| < \infty$.

Furthermore, $Ee^{\nu_0 X_1} = 1$ implies (by [32], Theorem 25.17) that $\int_{|x|>1} e^{\nu_0 x} \times \Pi(dx) < \infty$ and thus (by differentiation) $-\Psi(-i\nu)$ is finite and strictly convex for $\nu < \nu_0$. This rules out the possible existence of an $\alpha > 0$ with $Ee^{\alpha X_1} < 1$ and $\overline{\Pi}_X^+ \in \mathcal{S}^{(\alpha)}$, because the latter implies $Ee^{(\alpha+\varepsilon)X_1} = \infty$ for all $\varepsilon > 0$, while the convexity of $-\Psi(-i\nu)$ means that $\alpha < \nu_0$. Thus the situation in [6] and ours are mutually exclusive.

PROPOSITION 5.3 (Relation between $\Pi_X$ and $\Pi_{\mathcal{H}}$, $\alpha > 0$). *Assume* (4.1). *Then* $\overline{\Pi}_X^+$ *belongs to* $\mathcal{L}^{(\alpha)}$ *for a given* $\alpha > 0$ *if and only if* $\overline{\Pi}_{\mathcal{H}}$ *does, in which case* $\overline{\Pi}_X^+(u) \sim \hat{\kappa}(-i\alpha)\overline{\Pi}_{\mathcal{H}}(u)$, *as* $u \to \infty$.

Define
$$A_-(x) = \overline{\Pi}_X^-(1) + \int_1^x \overline{\Pi}_X^-(y)\, dy, \qquad x \geq 1,$$

and let "$\asymp$" in a relationship denote that ratio of the two sides is bounded away from zero and infinity, over the indicated range of the variable.

PROPOSITION 5.4 (Relation between $\Pi_X$ and $\Pi_{\mathcal{H}}$, $\alpha = 0$). *Assume* (4.1) *and* $\overline{\Pi}_X^+ \in \mathcal{L}^{(0)}$.

(i) *If* $\int_1^\infty \overline{\Pi}_X^-(y)\, dy = \infty$, *then*

$$\overline{\Pi}_{\mathcal{H}}(u) \asymp \int_{(1,\infty)} \left(\frac{y}{A_-(y)}\right) \Pi_X(u+dy), \qquad u \to \infty. \tag{5.3}$$

(ii) *If* $\int_1^\infty \overline{\Pi}_X^-(y)\, dy < \infty$, *then*

$$\overline{\Pi}_{\mathcal{H}}(u) \asymp \int_{(u,\infty)} \overline{\Pi}_X^+(y)\, dy, \qquad u \to \infty. \tag{5.4}$$

REMARK 5.5. (i) By [11], $\lim_{t \to \infty} X_t = -\infty$ a.s. if and only if

$$\int_1^\infty \overline{\Pi}_X^-(y)\, dy = \infty \quad \text{and} \quad \int_{(1,\infty)} \left(\frac{y}{A_-(y)}\right) \Pi_X(dy) < \infty, \tag{5.5}$$
$$\text{or} \quad 0 < -EX_1 \leq E|X_1| < \infty.$$

Thus the integral on the right-hand side of (5.3) is finite under (4.1).

(ii) We can apply (4.9), (5.3) and (5.4) as follows. Denote the right-hand side of (5.3) or (5.4) by $\overline{B}_0(u)$, a finite, nonincreasing function on $(0, \infty)$. Suppose there are functions $a(u) \to \infty$ as $u \to \infty$ and $b(x) \to 0$ as $x \to \infty$ such that, for each $x > 0$,

$$\frac{\overline{B}_0(u + xa(u))}{\overline{B}_0(u)} \asymp b(x), \qquad u \to \infty. \tag{5.6}$$

Note that $a(u)$ and $b(x)$ are defined in terms of the Lévy characteristics of $X$, rather than of $\mathcal{H}$. Assume $\overline{\Pi}_X^+ \in \mathcal{L}^{(0)}$. Then by (5.3) or (5.4) and (4.9), and using the uniformity of convergence in (4.9), we have, for each $x > 0$,

$$P(X_{\tau(u)} - u > xa(u) \mid \tau(u) < \infty) \asymp b(x), \qquad u \to \infty. \tag{5.7}$$

This gives the asymptotic order of magnitude of the overshoot, when normalized by $a(u)$; it tells us that $(X_{\tau(u)} - u)/a(u)$ is tight as $u \to \infty$, conditional on $\tau(u) < \infty$. It is the counterpart of the corresponding result in (1.5) of [3], except that [3] obtains a limit rather than an order of magnitude estimate, as a result of its more restrictive (finite mean and maximum domain of attraction) but more informative assumptions. We can likewise strengthen (5.3), replacing "≍" by "∼," under stronger assumptions, using methods such as those of [34], for example. We omit further details of this here. When $\overline{\Pi}_X^+ \in \mathcal{L}^{(0)}$, it is shown in the proof of Proposition 5.4 that $\overline{\Pi}_{\mathcal{H}}(u)/\overline{\Pi}_X^+(u) \to \infty$ as $u \to \infty$, so we cannot replace the hypothesis $\overline{\Pi}_{\mathcal{H}} \in \mathcal{L}^{(0)}$ by $\overline{\Pi}_X^+ \in \mathcal{L}^{(0)}$ in Theorem 4.4.

(iii) Theorems 4.1, 4.2 and 4.4 generalize corresponding results of Asmussen and Klüppelberg [3] concerning the classical insurance risk process. In their case the limit d.f. $G$ reduces to a generalized Pareto distribution: for $\alpha = 0$, the normalizing function $a(u)$ from (5.7) is the well-known *auxiliary function* in extreme value theory (see [16], Chapter 3), and $G$ is a Pareto distribution; for $\alpha > 0$, the normalizing function degenerates to the constant $\alpha$, and $G$ is the standard exponential d.f.

**6. Examples.** In this section we shall consider examples, all of which have the feature that $X$ is spectrally positive; that is to say, we assume $\overline{\Pi}_X^- = 0$. This case is very tractable and allows us to derive quite explicit expressions which generalize well-known results in collective risk theory. It is the case of most direct interest in insurance applications. As before, $X$ also drifts to $-\infty$ a.s. For such processes, from Remark 2.11 we have that the downward ladder height process is simply a negative linear drift, $\hat{H}_t = -t$, $\overline{\Pi}_{\mathcal{H}}(u) = \int_u^\infty \overline{\Pi}_X^+(y)\,dy$ (finite), for $u > 0$, $E|X_1| < \infty$, and $EX_1 < 0$.

ASSUMPTION 6.1. Fix an $\alpha \geq 0$. When $\alpha > 0$, assume that $\overline{\Pi}_X^+ \in \mathcal{S}^{(\alpha)}$, and when $\alpha = 0$, assume that $\overline{\Pi}_{\mathcal{H}}$ belongs to $\mathcal{S}^{(0)}$. (Thus in either case we have $\overline{\Pi}_{\mathcal{H}} \in \mathcal{S}^{(\alpha)}$.)

Suppose that $X$ has Laplace exponent $\phi(\theta)$ for $\theta \in \mathbb{R}$ such that

$$E(e^{\theta X_t}) = e^{\phi(\theta)t}.$$

The introduction of $\phi$ conveniently connects with existing literature on one-sided Lévy processes. When $\phi$ is finite, $\phi$ and $\Psi$ are related through the identity $\phi(\theta) = -\Psi(-i\theta)$. Under Assumption 6.1, $\phi(\theta)$ is finite for $\theta \in (-\infty, \alpha]$ and infinite for $\theta \in (\alpha, \infty)$. Noting that $-X$ is spectrally negative, we can extract the following facts from [4], Chapter VII, [5] and [32]: the function $\phi(\theta)$ is strictly convex on $(-\infty, \alpha]$, passes though the origin, has $\lim_{\theta \to -\infty} \phi(\theta) = +\infty$, and the drift of $X$ is given by the left-hand derivative $\phi'(0-) = EX_1$, which is finite and strictly negative.

Using (2.9) and taking advantage of the fact that the downward ladder height process is simply a linear drift, we can identify $q$ as

$$q = \lim_{\theta \downarrow 0} \frac{\Psi(i\theta)}{\hat{\kappa}(i\theta)} = \lim_{\theta \downarrow 0} \frac{-\phi(-\theta)}{-\log E(e^{-\theta \hat{H}_1})} = \lim_{\theta \downarrow 0} \frac{\phi(-\theta)}{\theta}.$$

We thus deduce that $q = -\phi'(0-) = |EX_1| < \infty$. Next note from (5.2) that, when $\alpha > 0$,

$$q - \log \delta_\alpha(\mathcal{H}) = \frac{-\phi(\alpha)}{\alpha},$$

and hence the condition $e^{-q}\delta_\alpha(\mathcal{H}) < 1$ reduces to the requirement that $\phi(\alpha) < 0$. (Recall that when $\alpha = 0$, the requirement $e^{-q}\delta_\alpha(\mathcal{H}) < 1$ is automatically satisfied.)

We can now read off the following conclusions from (2.13) and Theorems 4.1 and 4.2.

THEOREM 6.2. *Suppose that $X$ is spectrally positive, drifts to $-\infty$ a.s., satisfies Assumption 6.1 for a given $\alpha \geq 0$, and has $\phi(\alpha) < 0$ if $\alpha > 0$. With the understanding that $-\phi(\alpha)/\alpha = -\phi'(0-) = |EX_1|$ when $\alpha = 0$, we have:*

(i) $P(\tau(u) < \infty) \sim |EX_1| \left(\dfrac{\alpha}{\phi(\alpha)}\right)^2 \displaystyle\int_u^\infty \overline{\Pi}_X^+(y)\, dy \qquad as\ u \to \infty;$

(ii) $\lim_{u \to \infty} P(X_{\tau(u)} - u > x | \tau(u) < \infty) = \overline{G}(x)$, *where*

$$\overline{G}(x) = \frac{e^{-\alpha x}}{|EX_1|}\left(\frac{-\phi(\alpha)}{\alpha} + \int_{(x,\infty)}(e^{\alpha y} - e^{\alpha x})\overline{\Pi}_X^+(y)\, dy\right);$$

(iii) $\displaystyle\lim_{u\to\infty} P(L_{\tau(u)} > t | \tau(u) < \infty) = e^{\phi(\alpha)t/\alpha}\left(1 - \frac{t\phi(\alpha)}{\alpha}\left(1 + \frac{\phi(\alpha)}{\alpha|EX_1|}\right)\right).$

Let us proceed to examine some specific spectrally positive models in more detail.

6.1. *Jump diffusion process.* Suppose Assumption 6.1 is in force and, further, that $X_t$ is of the form

(6.1) $$X_t = \sigma B_t + \sum_{i=1}^{N_t} Y_i - \gamma t, \qquad t \geq 0,$$

where $\gamma > 0$ and $\sigma > 0$ are constants, $B_t$ is a standard Brownian motion, $N_t$ is a Poisson process of rate $\lambda$, the $Y_i$ are a.s. positive i.i.d. r.v.'s with d.f. $F$ and all processes are independent. In the context of insurance risk theory, this process is called a *risk process perturbed by Brownian motion*; see [12] and [22]. $X_t$ can drift to $-\infty$ only if $EY_1 = \mu < \infty$, so assume this. For this process we have

$$(6.2) \quad \phi(\theta) = -\theta\gamma + \sigma^2\theta^2/2 + \lambda \int_{[0,\infty)} (e^{\theta x} - 1) F(dx),$$

and this is finite for $\theta \leq \alpha$ by Assumption 6.1. Also, $\overline{\Pi}_X^+(x) = \lambda \overline{F}(x)$, so

$$(6.3) \quad \int_u^\infty \overline{\Pi}_X^+(y)\,dy = \lambda \int_u^\infty \overline{F}(y)\,dy = \lambda\mu\overline{F}_I(u), \qquad u \geq 0,$$

where $F_I$ is the integrated tail d.f. as defined in (1.4).

(i) Take $\alpha > 0$. By Assumption 6.1, $\overline{\Pi}_X^+ \in \mathcal{S}^{(\alpha)}$ for the specified $\alpha$, so $F \in \mathcal{S}^{(\alpha)}$; thus we have $F \in \mathcal{L}^{(\alpha)}$, or, equivalently, $\overline{F} \circ \log$ is regularly varying with index $-\alpha$. Then by Karamata's theorem (see [8], page 28), we have $\lim_{u\to\infty} \overline{F}(u)/\int_u^\infty \overline{F}(y)\,dy = \alpha$. Hence by tail equivalence also $\overline{F}_I \in \mathcal{S}^{(\alpha)}$ and $\delta_\alpha(F_I) = (\delta_\alpha(F) - 1)/(\mu\alpha) < \infty$. It follows from (6.2) that

$$-\phi(\alpha) = \gamma\alpha - \sigma^2\alpha^2/2 - \lambda \int_{[0,\infty)} (e^{\alpha x} - 1) F(dx)$$
$$= \gamma\alpha - \sigma^2\alpha^2/2 - \lambda(\delta_\alpha(F) - 1)$$
$$= \gamma\alpha - \sigma^2\alpha^2/2 - \lambda\mu\alpha\delta_\alpha(F_I),$$

and this is positive if and only if (recall that $\rho = \mu\lambda/\gamma$)

$$(6.4) \quad \rho\delta_\alpha(F_I) + \frac{\sigma^2}{2}\frac{\alpha}{\gamma} < 1,$$

which we will assume to be the case. Note that this implies $\rho < 1$ and hence $\lim_{t\to\infty} X_t = -\infty$ a.s., because Wald's lemma and (6.1) show that $-EX_1 = \gamma - \lambda\mu = \gamma(1 - \rho) > 0$. Finally we have via substitution in Theorem 6.2(i) with the help of (6.3) and (6.4) that, as $u \to \infty$,

$$(6.5) \quad \begin{aligned} P(\tau(u) < \infty) &\sim \frac{(1-\rho)\rho}{(1 - \sigma^2\alpha/(2\gamma) - \rho\delta_\alpha(F_I))^2} \overline{F}_I(u) \\ &\sim \frac{(1-\rho)\rho}{\mu\alpha(1 - \sigma^2\alpha/(2\gamma) - \rho\delta_\alpha(F_I))^2} \overline{F}(u). \end{aligned}$$

This holds under Assumption 6.1 and (6.4). Similarly, we can obtain a quite explicit expression for the overshoot limit distribution from Theorem 6.2(ii), calculable once the (incomplete) moment generating function is calculated.

(ii) Take $\alpha = 0$. Our assumption is now that

$$\overline{\Pi}_{\mathcal{H}}(\cdot) = \lambda \int_\cdot^\infty \overline{F}(y)\,dy \in \mathcal{S}^{(0)}.$$

Assuming $\rho = \mu\lambda/\gamma < 1$ and again applying Theorem 6.2(i), we get

$$P(\tau(u) < \infty) \sim \frac{\rho}{1-\rho}\overline{F}_I(u), \qquad u \to \infty, \tag{6.6}$$

in which the effect of the Brownian component has washed out. Equation (6.6) is the same as for the classical case; see [18] and [16], Section 1.4. In this case $\alpha = 0$, (4.6) simply tells us that the overshoot above $u$ tends to $\infty$ in probability as $u \to \infty$, as we expect from the heavy-tailedness in the positive direction. To sharpen the result we use Theorem 4.4 and the arguments in Remark 5.5. With $\overline{B}_0(u) = \lambda \int_u^\infty \overline{F}(y)\,dy$, choose $a(u)$ to satisfy (5.6), so that (5.7) holds. This parallels the development of Asmussen and Klüppelberg [3] for essentially the same model (without the Brownian component). They require a maximum domain of attraction condition, which gives "$\to$" in (5.7), whereas our more general analysis only gives "$\asymp$."

REMARK 6.3. The last result can be viewed as a robustness result in the sense of suggesting how far we can move away from a specific model without changing the asymptotic ruin probability: to a random walk with subexponential claims, we can add a diffusion term without changing the ruin probability. This effect has been investigated in a more general framework by Embrechts and Samorodnitsky [17]. See also [29]. Our next example also has an interpretation in this sense.

6.2. *Stable process with jumps and drift.* In this example we suppose $X_t$ is of the form

$$X_t = S_t^{(p)} + \sum_{i=1}^{N_t} Y_i - \gamma t, \qquad t \geq 0, \tag{6.7}$$

where $\gamma > 0$, $S_t^{(p)}$ is a stable Lévy motion with index $p \in (1,2)$, and the variables $\{Y_i : i \geq 1\}$ are as before (thus, with $EY_1 = \mu < \infty$). It follows from [38] that

$$\phi(\theta) = -\gamma\theta + \int_{[0,\infty)} (e^{\theta x} - 1 - \theta x)\frac{p(p-1)}{\Gamma(2-p)x^{1+p}}\,dx + \lambda \int_{[0,\infty)} (e^{\theta x} - 1)F(dx)$$

$$= -\gamma\theta + (-\theta)^p + \lambda \int_{[0,\infty)} (e^{\theta x} - 1)F(dx), \qquad \theta \leq 0.$$

For this example the mgf $\phi(\theta)$ is finite only if $\theta \leq 0$, so we only consider the case $\alpha = 0$; that is, we assume $\overline{\Pi}_\mathcal{H} \in \mathcal{S}^{(0)}$. The process has no downward jumps ($\beta = 1$ in the notation of [4], Chapter VIII, page 217). This model has been considered by Furrer [21] and Schmidli [33] and is in the ruin context called a *risk process perturbed by p-stable Lévy motion*.

Again assume $\rho = \mu\lambda/\gamma < 1$. By differentiation, $q = -\phi'(0-) = \gamma - \lambda\mu$. The Lévy measure of $X$ satisfies

$$\overline{\Pi}_X^+(x) = \frac{(p-1)}{\Gamma(2-p)} x^{-p} + \lambda \overline{F}(x) \quad \text{and} \quad \overline{\Pi}_X^-(x) = 0, \qquad x > 0,$$

and, further,

$$\overline{\Pi}_{\mathcal{H}}(u) = \int_u^\infty \overline{\Pi}_X^+(x)\,dx = \frac{u^{-(p-1)}}{\Gamma(2-p)} + \lambda \int_u^\infty \overline{F}(x)\,dx, \qquad u > 0.$$

We distinguish three different cases: suppose

$$\lim_{x \to \infty} x^p \overline{F}(x) = \begin{cases} 0, \\ c \in (0, \infty), \quad \text{or} \\ \infty. \end{cases}$$

From l'Hôpital we get, corresponding to the above cases,

$$\overline{\Pi}_{\mathcal{H}}(u) \sim \begin{cases} \dfrac{u^{-(p-1)}}{\Gamma(2-p)}, \\ \left(\dfrac{1}{\Gamma(2-p)} + \dfrac{\lambda c}{p-1}\right) u^{-(p-1)}, \quad \text{or} \\ \lambda\mu \overline{F}_I(u). \end{cases}$$

This means that we are in the same situation as in the classical subexponential case (i.e., when $\alpha = 0$), but have two different regimes depending on whether the tail of the claim size distribution is heavier or lighter than that of the stable perturbation.

Consequently, for $F$ with tails lighter than or similar to $x^{-p}$ (i.e., Cases 1 and 2), Theorem 6.2 gives

$$P(\tau(u) < \infty) \sim \frac{C}{(\gamma - \lambda\mu)} u^{-(p-1)}, \qquad u \to \infty,$$

with $C = 1/\Gamma(2-p)$ or $C = 1/\Gamma(2-p) + \lambda c/(p-1)$. If $F$ is heavier tailed than $x^{-p}$ (Case 3), we again get (6.6), with $\rho = \mu\lambda/\gamma$.

6.3. *Notes and comments.* All models considered in the insurance literature so far have entailed very specific Lévy processes; in particular, of course, the classical compound Poisson model as introduced in Section 1 has gained a lot of attention. In [3], page 106, it is suggested that "by a discrete skeleton argument," it may be possible to extend their random walk results to a general Lévy process. There are some difficulties with transferring results in this way, however, to do with relating the passage time above a level $u$ of the discrete process $\{X_n\}_{n=1,2,\ldots}$ to the continuous time version $\tau(u)$, or, more generally, relating the ladder processes and corresponding Lévy measures in a useful way. An alternative approach via a path decomposition of the Lévy process into drift, diffusion, "small jump" and "large jump" processes seems to run into similar problems. Our direct approach to the ladder properties of the Lévy process itself, with the help of the Bertoin and Doney [6] and Vigon [36] techniques, avoids these considerations and provides a basis for further developments of the theory.

In previous investigations, apart from estimates for the ruin probability, interest has mostly been concentrated on working out joint limiting distributions of the ruin time $\tau(u)$, the surplus $X_{\tau(u)-}$ before ruin, and/or the overshoot $X_{\tau(u)} - u$ after ruin. This problem was also considered for the classical model perturbed by a Brownian motion as in Section 6.1; see [23]. A more general approach than this is pursued by Huzak, Perman, Arvoje and Vondracek [26]. They consider a perturbed risk process defined as

$$X_t = Z_t + C_t - \gamma t, \qquad t \geq 0,$$

where $C$ is a subordinator representing the total claim amount process; it has only upward jumps. The perturbation $Z$ is a Lévy process, independent of $C$, which is also spectrally positive with zero expectation. In this model the Pollacek–Khinchine formula can again be given quite explicitly, stating for the survival probability

$$1 - \psi(u) = (1 - \rho) \sum_{n=1}^{\infty} \rho^n (G^{*(n+1)} * H^{*n})(u), \qquad u \geq 0.$$

The parameter $\rho$ can again be specified, and $G$ and $H$ are d.f.'s, where $G$ can be identified as the d.f. of the absolute supremum of the process $\{Z_t - \gamma t : t \geq 0\}$ and $H$ is the integrated tail d.f. of the jumps of $C$. The main concern in [26] is to analyze the supremum and ladder height processes of $X$, $Z$ and $C$. Since $X$ is a Lévy process, our results also apply to it, and analysis along our lines can be carried out; but we do not proceed further here.

Finally we remark that all of our previous general results have exact random walk analogues too, assuming only that the random walk drifts to $-\infty$ a.s., and that the distribution of the increments satisfies similar subexponential/convolution equivalence conditions and a non-Cramér condition as we imposed for the Lévy process. The results can even be strengthened slightly. Since the proofs for the discrete time case use the same ideas, and are even a little simpler, we omit the details.

**7. Proofs.** We need a couple of technical lemmas. The first is a minor modification of some working out in [14].

LEMMA 7.1. *Let $\alpha \geq 0$ and let the d.f. $\nu(\cdot) \in \mathcal{S}^{(\alpha)}$. Then, for each $x \geq 0$,*

$$\lim_{a \to \infty} \limsup_{u \to \infty} \int_{(a, u+x-a]} \frac{\overline{\nu}(u+x-y)}{\overline{\nu}(u)} \nu(dy) = 0. \tag{7.1}$$

*Further, the convergence in* (7.1) *is uniform in* $x \geq 0$.

PROOF. Write
$$\overline{\nu^{*2}}(z) = \overline{\nu}(z) + \overline{\nu} * \nu(z) = \overline{\nu}(z) + \int_{(0,z]} \overline{\nu}(z-y)\nu(dy), \qquad z \geq 0.$$

For $a > 0$ and $z > 2a$, split up the convolution integral into integrals over $(0, a], (a, z - a), [z - a, z)$ and use partial integration on the last integral. This gives the identity

$$(7.2) \qquad \overline{\nu^{*2}}(z) = \left(2\int_{(0,a]} + \int_{(a,z-a)}\right)\overline{\nu}(z-y)\nu(dy) + \overline{\nu}(a)\overline{\nu}(z-a),$$

from which we see that

$$(7.3) \quad \begin{aligned} &\int_{(a,u+x-a)} \overline{\nu}(u+x-y)\nu(dy) \\ &\qquad = \overline{\nu^{*2}}(u+x) - 2\int_{(0,a]} \overline{\nu}(u+x-y)\nu(dy) - \overline{\nu}(a)\overline{\nu}(u+x-a). \end{aligned}$$

Now $\nu(\cdot) \in \mathcal{S}^{(\alpha)}$, so $\overline{\nu^{*2}}(u)/\overline{\nu}(u) \to 2\delta_\alpha(\nu)$ and $\overline{\nu}(u-z) \sim e^{\alpha z}\overline{\nu}(u)$ as $u \to \infty$. Divide by $\overline{\nu}(u)$ in (7.3) and let $u \to \infty$, using dominated convergence, to get the limit

$$2e^{-\alpha x}\int_{(0,\infty)} e^{\alpha y}\nu(dy) - 2\int_{(0,a]} e^{\alpha(y-x)}\nu(dy) - e^{\alpha(a-x)}\overline{\nu}(a)$$
$$= 2e^{-\alpha x}\int_{(a,\infty)} e^{\alpha y}\nu(dy) - e^{-\alpha x}e^{\alpha a}\overline{\nu}(a).$$

The finiteness of $\delta_\alpha(\nu)$ implies $\lim_{a \to \infty} e^{\alpha a}\overline{\nu}(a) = 0$, so the last expression tends to 0 as $a \to \infty$.

For the uniformity, note that (7.1) with $x = 0$ gives $\int_{(a,u-a]} \overline{\nu}(u-y)\nu(dy) \leq \varepsilon\overline{\nu}(u)$ once $a \geq a_0(\varepsilon)$ and $u \geq u_0(a)$. Then if $x \geq 0$ and $u_x = u+x$, $\int_{(a,u_x-a]} \overline{\nu}(u_x - y)\nu(dy) \leq \varepsilon\overline{\nu}(u)$ if $a \geq a_0(\varepsilon)$ and $u_x \geq u_0(a)$, certainly if $a \geq a_0(\varepsilon)$ and $u \geq u_0(a)$. $\square$

We shall use the nonlattice part of the following lemma; for completeness, we also include the lattice case. It is simply a re-presentation of the defining properties for $\mathcal{L}^{(\alpha)}$, and we omit the proof.

LEMMA 7.2. *For* $\alpha > 0$, $G \in \mathcal{L}^{(\alpha)}$ *is equivalent to*

$$\lim_{u \to \infty} \frac{G(u, u+h)}{G(u, u+1)} = \frac{1 - e^{-\alpha h}}{1 - e^{-\alpha}},$$

*where the limit is through values $u$ in $\mathbb{R}$ or $\mathbb{N}$, and for all $h > 0$ or for all $h \in \mathbb{N}$, for the nonlattice and lattice case, respectively. This in turn is*

*equivalent to saying that $G(u+dy)/\overline{G}(u)$ converges weakly to an exponential distribution with parameter $\alpha$ or to a geometric distribution with parameter $e^{-\alpha}$, respectively.*

PROOF OF THEOREM 4.1. Fix $\alpha \geq 0$ and assume (4.1)–(4.3). The forward part of the theorem follows from (3.8), together with the use of (3.7).

For the converse part, assume (4.1), let $\overline{C}(u) = P(\tau(u) < \infty) = q\overline{V}(u)$, $u > 0$, so that $\overline{C}(u)$ is the tail of a d.f. $C$, and suppose $C \in \mathcal{S}^{(\alpha)}$. Thus $\delta_\alpha(C) = q\delta_\alpha(V) < \infty$, and since

$$\begin{aligned}(7.4)\quad \delta_\alpha(V) &= \int_{[0,\infty)} e^{\alpha y} \int_0^\infty e^{-qt}\,dt \cdot P(\mathcal{H}_t \in dy) \\ &= \int_0^\infty e^{-qt}\,dt \cdot E(e^{\alpha \mathcal{H}_t}) \\ &= \int_0^\infty e^{-qt}\,dt \cdot \delta_\alpha^t(\mathcal{H}),\end{aligned}$$

(4.3) holds. Now $\overline{C}(u)$ satisfies

$$\overline{C}(u) = \overline{C}_q(u) := q\int_0^\infty e^{-qt} P(\mathcal{H}_t > u)\,dt = P(\mathcal{H}_{\mathbf{e}_q} > u), \qquad u > 0,$$

where $\mathbf{e}_q$ is an independent exponential variable with parameter $q$, and, since $C_q \in \mathcal{S}^{(\alpha)}$, we have by (3.4), for each $k \in \mathbb{N}$,

$$(7.5)\quad \lim_{u \to \infty} \frac{\overline{C_q^{*k}}(u)}{\overline{C_q}(u)} = k\delta_\alpha^{k-1}(C_q).$$

Using the fact that $\mathcal{H}$ has stationary independent increments, we have that $\overline{C_q^{*k}}(u) = P(\mathcal{H}_{\mathbf{e}_q^k} > u)$, where $\mathbf{e}_q^k$ is the sum of $k$ independent exponential r.v.'s, each with parameter $q$. That is to say,

$$\overline{C_q^{*k}}(u) = \frac{q^k}{(k-1)!} \int_0^\infty t^{k-1} e^{-qt} P(\mathcal{H}_t > u)\,dt, \qquad u > 0.$$

Thus by (7.5),

$$\lim_{u \to \infty} \frac{q^k}{(k-1)!\overline{C}_q(u)} \int_0^\infty t^{k-1} e^{-qt} P(\mathcal{H}_t > u)\,dt = k\delta_\alpha^{k-1}(C_q).$$

Multiplying both sides of this by $(1-\lambda/q)^{k-1}$, with $q(1-1/\delta_\alpha(C_q)) < \lambda < q(1+1/\delta_\alpha(C_q))$ [so that $|1-\lambda/q| < 1/\delta_\alpha(C_q)$], and summing over $k \in \mathbb{N}$, gives

$$\begin{aligned}(7.6)\quad \lim_{u\to\infty} \frac{1}{\overline{C}_q(u)} &\int_0^\infty e^{-\lambda t} P(\mathcal{H}_t > u)\,dt \\ &= \frac{1/q}{(1-(1-\lambda/q)\delta_\alpha(C_q))^2} \\ &= \frac{(q - \log \delta_\alpha(\mathcal{H}))^2}{q(\lambda - \log \delta_\alpha(\mathcal{H}))^2},\end{aligned}$$

because $\delta_\alpha(C_q) = q\delta_\alpha(V) = q/(q - \log \delta_\alpha(\mathcal{H}))$. Relation (7.6) is valid for $q(1 - 1/\delta_\alpha(C_q)) < \lambda < q(1 + 1/\delta_\alpha(C_q))$. It means that

$$\overline{C}_\lambda(u) = P(\mathcal{H}_{\mathbf{e}_\lambda} > u) = \lambda \int_0^\infty e^{-\lambda t} P(\mathcal{H}_t > u)\, dt \sim c\overline{C}_q(u),$$

for some $c > 0$ and hence $C_\lambda$ is in $\mathcal{S}^{(\alpha)}$ for $\lambda$ in the indicated range. So by repeating the above argument with $q$ replaced by a $\lambda_0 \in (q, q(1 + 1/\delta_\alpha(C_q)))$ [for which one should note that $\delta_\alpha(C_{\lambda_0}) < q^{-1}\lambda_0\delta_\alpha(C_q)$ and hence that $\lambda_0(1 + 1/\delta_\alpha(C_{\lambda_0})) > q(1 + 1/\delta_\alpha(C_q))$], we can extend the upper limit of the range of applicability of (7.6). Continuing in this way, we see that (7.6) holds for all $\lambda > q(1 - 1/\delta_\alpha(C_q))$. So we can write, for all large $\lambda$,

$$\lim_{u\to\infty} \frac{1}{\overline{C}_q(u)} \int_0^\infty e^{-\lambda t} P(\mathcal{H}_t > u)\, dt$$
$$= \frac{(q - \log \delta_\alpha(\mathcal{H}))^2}{q} \int_0^\infty t e^{-(\lambda - \log \delta_\alpha(\mathcal{H}))t}\, dt.$$

Then by the continuity theorem for Laplace transforms ([20], page 433), we get

$$\lim_{u\to\infty} \frac{P(\mathcal{H}_t > u)}{\overline{C}_q(u)} = \frac{(q - \log \delta_\alpha(\mathcal{H}))^2}{q} t\delta_\alpha^t(\mathcal{H}).$$

By tail equivalence this means that $P(\mathcal{H}_1 > u) \in \mathcal{S}^{(\alpha)}$. □

PROOF OF THEOREM 4.2. Fix $\alpha \geq 0$, and suppose throughout that (4.1)–(4.3) hold.

(i) Take $x > 0$ and $a > 0$, choose $u > 2a$, and write, from Corollary 2.9(i),

$$(7.7) \quad P(X_{\tau(u)} - u > x, \tau(u) < \infty) = \left(\int_{(0,a]} + \int_{(a,u]}\right) \overline{\Pi}_{\mathcal{H}}(u + x - y) V(dy)$$
$$=: A_u + B_u.$$

By (3.1) we have

$$\lim_{u\to\infty} \frac{\overline{\Pi}_{\mathcal{H}}(u - y)}{\overline{\Pi}_{\mathcal{H}}(u)} = e^{\alpha y}, \qquad y \in \mathbb{R}.$$

In $A_u$, $y \leq a$, so the integrand is dominated by

$$\overline{\Pi}_{\mathcal{H}}(u + x - a) \leq \overline{\Pi}_{\mathcal{H}}(u - a) \leq 2e^{\alpha a}\overline{\Pi}_{\mathcal{H}}(u), \qquad u \geq u_0(a),$$

for some $u_0(a)$ large enough. Thus by dominated convergence,

$$\lim_{u\to\infty} \frac{A_u}{\overline{\Pi}_{\mathcal{H}}(u)} = \int_{(0,a]} e^{\alpha(y-x)} V(dy),$$

and as the convergence of monotone functions, the convergence is uniform in $x \geq 0$. So by (3.8),

$$\lim_{u \to \infty} \frac{A_u}{P(\tau(u) < \infty)} = \frac{1}{q\delta_\alpha^2(V)} \int_{(0,a]} e^{\alpha(y-x)} V(dy)$$

$$= \frac{e^{-\alpha x}}{q\delta_\alpha^2(V)} \Big(\delta_\alpha(V) - \int_{(a,\infty)} e^{\alpha y} V(dy)\Big).$$

Since $\delta_\alpha(V) < \infty$, as shown in the proof of Theorem 4.1, when $a \to \infty$ the integral here tends to 0 and, with (7.4), we get the first two terms in (4.6).

Next we deal with $B_u$, in (7.7). Integration by parts gives

$$
\begin{aligned}
B_u &= \overline{\Pi}_{\mathcal{H}}(u+x-a)\overline{V}(a) - \overline{\Pi}_{\mathcal{H}}(x)\overline{V}(u) \\
&\quad + \int_{[x,u+x-a)} \overline{V}(u+x-y)\Pi_{\mathcal{H}}(dy) \\
(7.8) \qquad &= \overline{\Pi}_{\mathcal{H}}(u+x-a)(\overline{V}(a) - \overline{V}(u)) \\
&\quad + \int_{[x,u+x-a)} (\overline{V}(u+x-y) - \overline{V}(u))\Pi_{\mathcal{H}}(dy) \\
&=: \overline{\Pi}_{\mathcal{H}}(u+x-a)(\overline{V}(a) - \overline{V}(u)) + C_u.
\end{aligned}
$$

When divided by $\overline{\Pi}_{\mathcal{H}}(u)$, the first term is dominated by $\overline{\Pi}_{\mathcal{H}}(u-a)\overline{V}(a)/\overline{\Pi}_{\mathcal{H}}(u)$, which tends to $e^{\alpha a}\overline{V}(a)$ as $u \to \infty$, and since $\delta_\alpha(V) < \infty$, we have $e^{\alpha a}\overline{V}(a) \to 0$ as $a \to \infty$.

Take $a > x > 0$ and $u + x > a$ and write

$$
(7.9) \qquad \begin{aligned}
\frac{C_u}{\overline{V}(u)} &= \Big(\int_{[x,a]} + \int_{(a,u+x-a)}\Big) \Big(\frac{\overline{V}(u+x-y)}{\overline{V}(u)} - 1\Big)\Pi_{\mathcal{H}}(dy) \\
&=: \frac{D_u}{\overline{V}(u)} + \frac{E_u}{\overline{V}(u)}.
\end{aligned}
$$

In the first term, the integrand is dominated by

$$\frac{\overline{V}(u-a)}{\overline{V}(u)} \leq 2e^{\alpha a}, \qquad u \geq u_1(a),$$

for $u_1(a)$ large enough, and a constant is integrable with respect to $\Pi_{\mathcal{H}}(dy)$ over $y \in (x, \infty)$, $x > 0$. Thus, by Proposition 2.5,

$$
\begin{aligned}
\frac{D_u}{P(\tau(u) < \infty)} &= \frac{D_u}{q\overline{V}(u)} \\
&\to \frac{1}{q}\int_{[x,a]} (e^{\alpha(y-x)} - 1)\Pi_{\mathcal{H}}(dy), \qquad u \to \infty,
\end{aligned}
$$

for each $a > 0$. This convergence of monotone functions is uniform in $x \in [\eta, \infty)$ for each $\eta > 0$. As $a \to \infty$, we get the last term on the right-hand side of (4.6).

The second term on the right-hand side of (7.9) is not negative, and since $\overline{V}(u) \leq c_0\, \overline{\Pi}_{\mathcal{H}}(u)$ for $u \geq u_2$, $u_2$ large enough, and some $c_0 > 0$, by (3.8), $E_u$ is bounded above by a constant multiple of

$$\text{(7.10)} \qquad \int_{(a,u+x-a)} \overline{\Pi}_{\mathcal{H}}(u+x-y)\Pi_{\mathcal{H}}(dy)$$

once $u + x - y \geq u_2$, and this is the case when $y < u + x - a$ if we choose $a > u_2$. Now since $\Pi_{\mathcal{H}}(\cdot) \neq 0$, we can choose $z_0 > 0$ such that $\overline{\Pi}_{\mathcal{H}}(z_0) > 0$. Also keep $a > z_0$. Then define

$$\nu(z) = \left(1 - \frac{\overline{\Pi}_{\mathcal{H}}(z)}{\overline{\Pi}_{\mathcal{H}}(z_0)}\right)\mathbf{1}_{\{z \geq z_0\}},$$

which is a (proper) d.f. on $[0, \infty)$ with tail $\overline{\nu}(z) = \overline{\Pi}_{\mathcal{H}}(z)/\overline{\Pi}_{\mathcal{H}}(z_0)$, $z > z_0$. The integral in (7.10) is, apart from a constant multiple,

$$\text{(7.11)} \qquad \int_{(a,u+x-a)} \overline{\nu}(u+x-y)\nu(dy).$$

The proof of part (i) is now complete with Lemma 7.1.

To prove part (ii), use the strong Markov property at the stopping time $L_t^{-1}$ to deduce

$$\begin{aligned}
P(L_{\tau(u)} > t, \tau(u) < \infty) &= P(H_t < u, \tau(u) < \infty) \\
&= E(\mathbf{1}_{(H_t < u)}P(\tau(u) < \infty|\mathcal{F}_{L_t^{-1}})) \\
&= E(\mathbf{1}_{(H_t < u)}P_{H_t}(\tau(u) < \infty)) \\
&= E(\mathbf{1}_{(H_t < u)}P(\tau(u - H_t) < \infty)) \\
&= E(\mathbf{1}_{(H_t < u)}P(\tau(u - H_t) < \infty); t < L_\infty) \\
&= e^{-qt}E(\mathbf{1}_{(\mathcal{H}_t < u)}P(\tau(u - \mathcal{H}_t) < \infty)) \\
&= qe^{-qt}\int_{(0,u)} \overline{V}(u-y)P(\mathcal{H}_t \in dy).
\end{aligned}$$

Write the last expression as

$$\text{(7.12)} \qquad qe^{-qt}\left(\int_{(0,a]} + \int_{(a,u-a]} + \int_{(u-a,u)}\right)\overline{V}(u-y)P(\mathcal{H}_t \in dy),$$

where $u > 2a > 0$, then divide it by $P(\tau(u) < \infty) = q\overline{V}(u)$ and let $u \to \infty$. By dominated convergence, the first term tends to

$$e^{-qt}\int_{(0,a]} e^{\alpha y}P(\mathcal{H}_t \in dy) = e^{-qt}\left(\delta_\alpha(\mathcal{H}_t) - \int_{(a,\infty)} e^{\alpha y}P(\mathcal{H}_t \in dy)\right),$$

and this tends to $e^{-qt}\delta_\alpha^t(\mathcal{H})$ as $a \to \infty$. (Recall that $t$ is kept constant in this proof.)

By Lemma 3.5, we can choose $a$ such that
$$\overline{V}(y) \leq c_1 \overline{\Pi}_{\mathcal{H}}(y) \leq c_1 P(\mathcal{H}_t > y)$$
for $y \geq a$, and some $c_1 > 0$, so the second integral in (7.12) is not larger than
$$c_1 \int_{(a, u-a]} P(\mathcal{H}_t > u - y) P(\mathcal{H}_t \in dy)$$
and, after division by $P(\mathcal{H}_t > u)$, this tends to 0 as $u \to \infty$ then $a \to \infty$ by Lemma 7.1.

Finally,
$$\begin{aligned}
(7.13) \quad & q \int_{(u-a, u)} \overline{V}(u - y) P(\mathcal{H}_t \in dy) \\
&= q\overline{V}(a) P(\mathcal{H}_t > u - a) - q\overline{V}(0) P(\mathcal{H}_t > u) \\
&\quad + q \int_{(0, a)} P(\mathcal{H}_t > u - y) V(dy),
\end{aligned}$$

while by (3.7) and (4.4), as $u \to \infty$,
$$\begin{aligned}
\frac{P(\mathcal{H}_t > u - y)}{\overline{V}(u)} &= \frac{P(\mathcal{H}_t > u - y)}{P(\mathcal{H}_1 > u - y)} \frac{P(\mathcal{H}_1 > u - y)}{P(\mathcal{H}_1 > u)} \frac{P(\mathcal{H}_1 > u)}{\overline{V}(u)} \\
&\to t \delta_\alpha^t(\mathcal{H}) e^{\alpha y} (q - \log \delta_\alpha(\mathcal{H}))^2 \\
&= c_t e^{\alpha y},
\end{aligned}$$

say. Thus the right-hand side of (7.13), when divided by $q\overline{V}(u)$, tends as $u \to \infty$ then $a \to \infty$ to
$$c_t(\delta_\alpha(V) - 1/q) = \frac{c_t \log \delta_\alpha(\mathcal{H})}{q(q - \log \delta_\alpha(\mathcal{H}))} = t\delta_\alpha^t(\mathcal{H})(q - \log \delta_\alpha(\mathcal{H})) \log \delta_\alpha(\mathcal{H})/q.$$

Thus the limit is
$$e^{-qt}(\delta_\alpha^t(\mathcal{H}) + t\delta_\alpha^t(\mathcal{H})(q - \log \delta_\alpha(\mathcal{H})) \log \delta_\alpha(\mathcal{H})/q),$$
which is the right-hand side of (4.7).

For part (iii), simply use Corollary 2.9 to write
$$P\left(X_{L_{L_{\tau(u)}^{-1}}^{-1}} \geq \phi, \tau(u) < \infty\right) = P(\tau(u) < \infty) - \int_{[0, \phi)} V(dy) \overline{\Pi}_{\mathcal{H}}(u - y),$$
divide by $P(\tau(u) < \infty)$, and take the limit as $u$ tends to infinity, using (4.4), to get (4.8). $\square$

PROOF OF THEOREM 4.4. Suppose that $\overline{\Pi}_{\mathcal{H}} \in \mathcal{S}^{(0)}$. Take $a = 0$ in (7.7) and (7.8) to see that
$$P(X_{\tau(u)} - u > x \mid \tau(u) < \infty) - \frac{\overline{V}(0)\overline{\Pi}_{\mathcal{H}}(u + x)}{P(\tau(u) < \infty)}$$

is bounded in modulus by

$$(7.14) \quad \frac{\overline{\Pi}_{\mathcal{H}}(u+x)\overline{V}(u)}{P(\tau(u)<\infty)} + \int_{[x,u+x)} \left( \frac{\overline{V}(u+x-y)-\overline{V}(u)}{P(\tau(u)<\infty)} \right) \Pi_{\mathcal{H}}(dy).$$

Since $\overline{V}(0) = 1/q$ and $P(\tau(u) < \infty) \sim \overline{\Pi}_{\mathcal{H}}(u)/q$ by (4.4), the first term in (7.14) converges to 0 (uniformly in $x \geq 0$) as $u \to \infty$, and it suffices to show that the integral converges to 0 as $u \to \infty$, where, in the denominator, we can replace $P(\tau(u) < \infty)$ by $\overline{\Pi}_{\mathcal{H}}(u)$ or by $\overline{V}(u)$. Take $a > 0$ and $u > a$ and write the integral in (7.14) as

$$\left( \int_{[x,u+x-a)} + \int_{[u+x-a,u+x)} \right) (\overline{V}(u+x-y) - \overline{V}(u)) \Pi_{\mathcal{H}}(dy).$$

The first integral on the right-hand side is the same one we dealt with in (7.8), called $C_u$, and consequently when divided by $\overline{V}(u)$ has the same limit as $C_u/\overline{V}(u)$ has, but with $\alpha = 0$, namely, 0. As observed there, the convergence is uniform in $x \in [\eta, \infty)$, $\eta > 0$.

Finally,

$$\int_{[u+x-a,u+x)} \frac{\overline{V}(u+x-y)}{\overline{\Pi}_{\mathcal{H}}(u)} \Pi_{\mathcal{H}}(dy)$$

$$\leq \overline{V}(0) \left( \frac{\overline{\Pi}_{\mathcal{H}}((u+x-a)-) - \overline{\Pi}_{\mathcal{H}}(u+x)}{\overline{\Pi}_{\mathcal{H}}(u)} \right)$$

and this tends to 0 as $u \to \infty$, for each $a > 0$, uniformly in $x \geq 0$. □

PROOF OF PROPOSITION 4.5. From Corollary 2.9(iii) and the remark in the proof thereof concerning an expression for $P(X_{\tau(u)} = u)$, we can write

$$P(L_{L_{\tau(u)^-}}^{-1} > \psi, \tau(u) < \infty) \leq \int_{(0,u)} V(dy;\psi) \overline{\Pi}_{\mathcal{H}}(u-y) + cV'(u)$$

$$= \left\{ \int_{(0,a)} + \int_{[a,u)} \right\} V(dy;\psi) \overline{\Pi}_{\mathcal{H}}(u-y) + cV'(u),$$

where $u > 2a > 0$. Also $P(\tau(u) < \infty) \sim c_3 \overline{\Pi}_{\mathcal{H}}(u)$ for some $c_3 > 0$, so it suffices for our purposes to divide by $\overline{\Pi}_{\mathcal{H}}(u)$. For the first integral,

$$\lim_{u \to \infty} \int_{(0,a)} V(dy;\psi) \frac{\overline{\Pi}_{\mathcal{H}}(u-y)}{\overline{\Pi}_{\mathcal{H}}(u)} \leq \lim_{u \to \infty} \frac{\overline{\Pi}_{\mathcal{H}}(u-a)}{\overline{\Pi}_{\mathcal{H}}(u)} V(a;\psi) \leq V(a;\psi).$$

However,

$$V(a;\psi) = \int_0^\infty dt \cdot e^{-qt} P(\mathcal{H}_t \leq a, \mathcal{L}_t^{-1} > \psi)$$

$$\leq \int_0^\infty dt \cdot e^{-qt} P(\mathcal{L}_t^{-1} > \psi)$$

$$= \frac{1}{q} P(\mathcal{L}_{\mathbf{e}_q}^{-1} > \psi) \to 0$$

as $\psi \to \infty$ because $\mathcal{L}_{\mathbf{e}_q}^{-1} < \infty$ almost surely. For the remaining terms, note that

$$\int_{[a,u)} V(dy;\psi) \overline{\Pi}_{\mathcal{H}}(u-y) + cV'(u) \leq \int_{[a,u)} V(dy) \overline{\Pi}_{\mathcal{H}}(u-y) + cV'(u)$$
$$= P\Big(X_{L_{L_{\tau(u)}-}^{-1}} \geq a, \tau(u) < \infty\Big).$$

Divide by $P(\tau(u) < \infty)$, take the limsup as $u$ tends to infinity, and then let $a \to \infty$. The result is zero by (4.8). □

PROOF OF PROPOSITION 5.1. Fix $\nu > 0$ and assume (4.1). Put $\theta = -i\nu$ in the Wiener–Hopf factorization (2.5) to get

(7.15) $\qquad -\log E(e^{\nu X_1}) = \log(e^{-q} E(e^{\nu \mathcal{H}_1})) \log E(e^{\nu \hat{H}_1})$

(when each side is finite). Now since $\lim_{t \to \infty} X_t = -\infty$ a.s., $E(e^{\nu \hat{H}_1})$ is always finite and less than 1 (recall our convention that $\hat{H}$ is negative), so $E(e^{\nu X_1})$ is finite and less than 1 if and only if $E(e^{-q} e^{\nu \mathcal{H}_1})$ is finite and less than 1. Thus the first equivalence in (5.1) holds, and the second equality in (5.2) holds. The second equivalence in (5.1) and the first equality in (5.2) follow from (7.4). Also, from (2.8), which is valid also for $\nu < 0$ when $\delta_\nu(\mathcal{H}) < \infty$, we get

$$\delta_\nu(\mathcal{H}) = E(e^{\nu \mathcal{H}_1}) = \exp\Big(\nu c + \int_{[0,\infty)} (e^{\nu y} - 1) \Pi_{\mathcal{H}}(dy)\Big),$$

giving the third equivalence in (5.1). The fourth equivalence in (5.1) follows from $Ee^{\nu X_1} < 1$ and (2.1). □

PROOF OF PROPOSITION 5.3. Fix $\alpha > 0$. Suppose first that $\overline{\Pi}_X^+ \in \mathcal{L}^{(\alpha)}$. (This part is based on the analogous version for random walks which appears in [7].) Using Theorem 2.10(ii), we have that

$$\Pi_{\mathcal{H}}(u, u+h) = \int_{(-\infty,0)} \Pi_X^+(u-y, u+h-y) \, d\hat{V}(y).$$

It is not difficult to justify integrating by parts to get

$$\frac{\Pi_{\mathcal{H}}(u, u+h)}{\overline{\Pi}_X^+(u)} = \int_{(0,\infty)} [\hat{V}(-y) - \hat{V}(-(y-h)^+)] \frac{\Pi_X(u+dy)}{\overline{\Pi}_X^+(u)}.$$

Since $\overline{\Pi}_X^+ \in \mathcal{L}^{(\alpha)}$, it follows that $\Pi_X(u+dy)/\overline{\Pi}_X^+(u)$ converges for $u \to \infty$ weakly to the exponential distribution. This, together with the fact that $\hat{V}$ is a renewal measure and hence the integrand in the last equality is uniformly bounded, implies that

$$\lim_{u\to\infty} \frac{\Pi_{\mathcal{H}}(u, u+h)}{\overline{\Pi}_X^+(u)} = \alpha(1 - e^{-\alpha h}) \int_{(0,\infty)} \hat{V}(-y) e^{-\alpha y}\, dy, \qquad h > 0.$$

Since the right-hand side is nonzero (recall Remark 2.11), Lemma 7.2 suffices to conclude that $\overline{\Pi}_{\mathcal{H}} \in \mathcal{L}^{(\alpha)}$.

Conversely, let $\overline{\Pi}_{\mathcal{H}} \in \mathcal{L}^{(\alpha)}$. Write Theorem 2.10(i) as

$$(7.16) \qquad \overline{\Pi}_X^+(u) = \int_{(u,\infty)} \overline{\Pi}_{\hat{H}}(u-y)\, d(\overline{\Pi}_{\mathcal{H}}(u) - \overline{\Pi}_{\mathcal{H}}(y)) + \hat{c}\Pi'_{\mathcal{H}}(u),$$

where the derivative is only present if $\hat{c} > 0$. By Fubini's theorem,

$$\overline{\Pi}_X^+(u) = \int_{(-\infty,0)} \Pi_{\mathcal{H}}(u, u-y) \Pi_{\hat{H}}(dy) + \hat{c}\Pi'_{\mathcal{H}}(u), \qquad u > 0.$$

Take $0 < h_1 < h_2$ and integrate both sides of the last equation to get

$$(7.17) \qquad \begin{aligned}\int_{u+h_1}^{u+h_2} \overline{\Pi}_X^+(z)\, dz &= \int_{(-\infty,0)} \left( \int_{h_1}^{h_2} \Pi_{\mathcal{H}}(u+z, u+z-y)\, dz \right) \Pi_{\hat{H}}(dy) \\ &\quad + \hat{c}(\overline{\Pi}_{\mathcal{H}}(u+h_1) - \overline{\Pi}_{\mathcal{H}}(u+h_2)).\end{aligned}$$

By dominated convergence,

$$\lim_{u\to\infty} \frac{1}{\overline{\Pi}_{\mathcal{H}}(u)} \int_{h_1}^{h_2} (\overline{\Pi}_{\mathcal{H}}(u+z) - \overline{\Pi}_{\mathcal{H}}(u+z-y))\, dz$$

$$= (1 - e^{\alpha y}) \int_{h_1}^{h_2} e^{-\alpha z}\, dz$$

$$= \alpha^{-1}(1 - e^{\alpha y})(e^{-\alpha h_1} - e^{-\alpha h_2}),$$

uniformly in $y \leq 0$. Thus, dividing (7.17) by $\overline{\Pi}_{\mathcal{H}}(u)$ and letting $u \to \infty$, we get

$$(7.18) \qquad \begin{aligned}\lim_{u\to\infty} \int_{h_1}^{h_2} &\left( \frac{\overline{\Pi}_X^+(u+z)}{\overline{\Pi}_{\mathcal{H}}(u)} \right) dz \\ &= \alpha^{-1}(e^{-\alpha h_1} - e^{-\alpha h_2})\left( \int_{(-\infty,0)} (1 - e^{\alpha y}) \Pi_{\hat{H}}(dy) + \hat{c}\alpha \right) \\ &= \alpha^{-1}(e^{-\alpha h_1} - e^{-\alpha h_2})\hat{\kappa}(-i\alpha).\end{aligned}$$

Finally, take any sequence $u''_n \to \infty$ and a subsequence $u'_n \to \infty$ such that, by Helly's theorem,

$$\lim_{u'_n\to\infty} \frac{\overline{\Pi}_X^+(u'_n + z)}{\overline{\Pi}_{\mathcal{H}}(u'_n)} = p(z), \qquad z > 0.$$

Using Fatou's lemma in (7.18) we see that the nonincreasing function $p(z)$ is finite for all $z > 0$, and then we can use dominated convergence in (7.18) to deduce that

$$\int_{h_1}^{h_2} p(z)\,dz = \alpha^{-1}(e^{-\alpha h_1} - e^{-\alpha h_2})\hat{\kappa}(-i\alpha).$$

Differentiating, we see that $p(z) = e^{-\alpha z}\hat{\kappa}(-i\alpha)$ for all $z > 0$, true also for all subsequences, and so

$$\overline{\Pi}_X^+(u+z) \sim e^{-\alpha z}\hat{\kappa}(-i\alpha)\overline{\Pi}_{\mathcal{H}}(u), \qquad u \to \infty.$$

Consequently, $\overline{\Pi}_X^+ \in \mathcal{L}^{(\alpha)}$, and the exact form of the asymptotic is also established.

□

PROOF OF PROPOSITION 5.4. By Remark 2.11, $\overline{\Pi}_{\hat{H}}(y) > 0$ for all $x > 0$. Next, Theorem 2.10(ii) gives

(7.19) $$\overline{\Pi}_{\mathcal{H}}(u) = \int_{(u,\infty)} \hat{V}(-(y-u))\Pi_X(dy), \qquad u > 0,$$

and [4], page 74, gives, for all $y \geq 0$,

(7.20) $$\hat{V}(-y) \asymp \frac{y}{\hat{c} + \hat{A}(y)} \asymp \frac{y}{\hat{A}(y)},$$

where

$$\hat{A}(x) = \int_0^x \overline{\Pi}_{\hat{H}}(y)\,dy$$

is nonzero for all $x > 0$. [The second asymptotic relation in (7.20) follows by considering cases $\hat{A}(\infty) = \infty$ or $\hat{A}(\infty) < \infty$.] Consequently, for all $u > 0$,

(7.21) $$\overline{\Pi}_{\mathcal{H}}(u) \asymp \int_{(u,\infty)} \left(\frac{y-u}{\hat{A}(y-u)}\right)\Pi_X(dy).$$

Now $\hat{A}(x)/x = \int_0^1 \overline{\Pi}_{\hat{H}}(xy)\,dy$ is nonincreasing, tends to 0 as $x \to \infty$ and has a positive (possibly infinite) limit as $x \to 0+$. Thus $a_0 := \lim_{x \to 0+}(x/\hat{A}(x))$ is finite (possibly 0).

Symmetrically to (7.19), we have

$$\overline{\Pi}_{\hat{H}}(u) = \int_{(u,\infty)} V(y-u)|\overline{\Pi}_X^-(dy)|, \qquad u > 0.$$

Now $V(y) \leq V(\infty) = 1/q$; thus $\overline{\Pi}_{\hat{H}}(u) \leq \overline{\Pi}_X^-(u)/q$ and it follows that

(7.22) $$\hat{A}(x) \leq \hat{A}(1) + \frac{1}{q}\int_1^x \overline{\Pi}_X^-(y)\,dy, \qquad x \geq 1.$$

If $\int_1^x \overline{\Pi}_X^-(y)\, dy = \infty$, then the right-hand side is asymptotic to $(1/q)\int_1^x \overline{\Pi}_X^-(y)\, dy$ as $x \to \infty$. For a reverse inequality in this case, choose $x_0 \geq 1$ so that $V(x_0) \geq 1/(2q)$. Then for $x > x_0$,

$$\hat{A}(x) = \hat{A}(1) + \int_1^x \overline{\Pi}_{\hat{H}}(z)\, dz = \hat{A}(1) + \int_1^x \int_{(z,\infty)} V(y-z)|\overline{\Pi}_X^-(dy)|\, dz$$

$$\geq \hat{A}(1) + \int_{x_0}^x \int_{(2z,\infty)} V(x_0)|\overline{\Pi}_X^-(dy)|\, dz$$

$$\geq \hat{A}(1) + \frac{1}{2q}\int_{x_0}^x \overline{\Pi}_X^-(2z)\, dz \asymp \int_1^x \overline{\Pi}_X^-(z)\, dz, \qquad x \to \infty.$$

Together these give

$$\hat{A}(x) \asymp \int_1^x \overline{\Pi}_X^-(y)\, dy \asymp A_-(x), \qquad x \to \infty,$$

thus $c_- A_-(x) \leq \hat{A}(x) \leq c_+ A_-(x)$ for some $0 < c_- \leq c_+ < \infty$ whenever $x \geq x_0 > 1$. Integration by parts in (7.21) gives

$$\overline{\Pi}_{\mathcal{H}}(u) \asymp a_0 \overline{\Pi}_X^+(u) + \int_0^\infty \overline{\Pi}_X^+(u+y)\, d\left(\frac{y}{\hat{A}(y)}\right), \qquad u \geq 0.$$

Now assume $\overline{\Pi}_X^+ \in \mathcal{L}^{(0)}$. Then Fatou's lemma applied to the last equation shows that $\overline{\Pi}_{\mathcal{H}}(u)/\overline{\Pi}_X^+(u) \to \infty$ as $u \to \infty$. Then

$$\int_{(u,u+x_0)} \left(\frac{y-u}{\hat{A}(y-u)}\right) \Pi_X(dy) \leq \left(\frac{x_0}{\hat{A}(x_0)}\right) \overline{\Pi}_X^+(u) = o(\overline{\Pi}_{\mathcal{H}}(u)), \qquad u \to \infty.$$

Thus, as $u \to \infty$,

$$\overline{\Pi}_{\mathcal{H}}(u) \asymp \int_{(u+x_0,\infty)} \left(\frac{y-u}{\hat{A}(y-u)}\right) \Pi_X(dy)$$

$$\asymp \int_{(u+x_0,\infty)} \left(\frac{y-u}{A_-(y-u)}\right) \Pi_X(dy)$$

$$\asymp \int_{(u+1,\infty)} \left(\frac{y-u}{A_-(y-u)}\right) \Pi_X(dy).$$

This proves (5.3) in case $\int_1^\infty \overline{\Pi}_X^-(y)\, dy = \infty$. If $\int_1^\infty \overline{\Pi}_X^-(y)\, dy < \infty$ then (7.22) gives $\hat{A}(\infty) < \infty$ and then (5.4) follows from (7.21). $\square$

**Acknowledgments.** We are grateful to Professor Ron Doney for helpful comments and discussions, and to a referee and an Associate Editor for suggestions which helped us to substantially improve the readability. We will take up the suggestion of the referee to present our results in a forthcoming survey paper aiming at more concrete actuarial questions.

C. KLÜPPELBERG  
CENTER FOR MATHEMATICAL SCIENCES  
MUNICH UNIVERSITY OF TECHNOLOGY  
D-80290 MUNICH  
GERMANY  
E-MAIL: cklu@ma.tum.de  
URL: www.ma.tum.de/stat/

A. E. KYPRIANOU  
DEPARTMENT OF MATHEMATICS  
UNIVERSITY OF UTRECHT  
P.O. BOX 80.010  
3500TA UTRECHT  
THE NETHERLANDS  
E-MAIL: kyprianou@math.uu.nl  
URL: www.math.uu.nl/people/kyprianou

R. A. Maller
Centre for Financial Mathematics, MSI
  and School of Finance & Applied Statistics
Australian National University
Canberra, ACT 0200
Australia
e-mail: ross.maller@anu.edu.au
url: www.ecel.uwa.edu.au/rmaller/Welcome.html